\documentclass[final,onefignum,onetabnum]{article}

\usepackage{array}
\usepackage{enumitem}
\usepackage{esvect}

\usepackage[caption=false]{subfig}
\usepackage{pgfplots}

\usepackage[colorlinks=true, bookmarksopen,
            pdfsubject={algorithms},
            linkcolor={blue},
            anchorcolor={black},
            citecolor={black},
            filecolor={magenta},
            menucolor={black},
            plainpages=false,pdfpagelabels,
            urlcolor={blue}]{hyperref}

\usepackage{algorithmic}

\usepackage{graphicx,epstopdf}

\usepackage{subfig}
\usepackage{amsmath}
\usepackage{amssymb}
\usepackage{graphicx}
\usepackage{dcolumn}
\usepackage{mathtools}
\usepackage{amsmath,amsfonts}
\usepackage{bm}
\usepackage{amsmath}
\usepackage{amssymb}
\usepackage{color}
\usepackage{float}
\usepackage{setspace}
\usepackage{tabularx}

\allowdisplaybreaks

\newcommand{\p}{\partial}
\newcommand{\f}{\frac}

\newcommand{\E}{ {\mathbb{E}} }
\newcommand{\ola}{\overleftarrow}

\newcommand{\be}{\begin{equation}}
\newcommand{\ee}{\end{equation}}

\def\ba{\begin{array}}                \def\ea{\end{array}}
\def\bel{\begin{equation}\label}      \def\ee{\end{equation}}

\colorlet{texcscolor}{blue!50!black}
\colorlet{texemcolor}{red!70!black}
\colorlet{texpreamble}{red!70!black}
\colorlet{codebackground}{black!25!white!25}

\date{}

\begin{document}

\title{A Unified Filter Method for Jointly Estimating State and Parameters of Stochastic Dynamical Systems via the Ensemble Score Filter}

\author{
Feng Bao\thanks{ Department of Mathematics, Florida State University, Tallahassee, Florida, \ ({\tt bao@math.fsu.edu}).} 
\and Zezhong Zhang \thanks{ Department of Mathematics, Florida State University, Tallahassee, Florida.}
\and Guannan Zhang  \thanks{ Computer Science and Mathematics Division, Oak Ridge National Laboratory, Oak Ridge, TN 37831, USA. 
 This manuscript has been authored by UT-Battelle, LLC, under contract DE-AC05-00OR22725 with the US Department of Energy (DOE). The US government retains and the publisher, by accepting the article for publication, acknowledges that the US government retains a nonexclusive, paid-up, irrevocable, worldwide license to publish or reproduce the published form of this manuscript, or allow others to do so, for US government purposes. DOE will provide public access to these results of federally sponsored research in accordance with the DOE Public Access Plan.}
      }

 \maketitle

\begin{abstract}
This paper tackles the intricate task of jointly estimating state and parameters in data assimilation for stochastic dynamical systems that are affected by noise and observed only partially.  While the concept of ``optimal filtering'' serves as the customary approach to estimate the state of the target dynamical system, traditional methods such as Kalman filters and particle filters encounter significant challenges when dealing with high-dimensional and nonlinear problems. When we also consider the scenario where the model parameters are unknown, the problem transforms into a joint state-parameter estimation endeavor. Presently, the leading-edge technique known as the Augmented Ensemble Kalman Filter (AugEnKF) addresses this issue by treating unknown parameters as additional state variables and employing the Ensemble Kalman Filter to estimate the augmented state-parameter vector. Despite its considerable progress, AugEnKF does exhibit certain limitations in terms of accuracy and stability. To address these challenges, we introduce an innovative approach, referred to as the United Filter. This method amalgamates a remarkably stable and efficient ensemble score filter (EnSF) for state estimation with a precise direct filter dedicated to online parameter estimation. Utilizing the EnSF's generative capabilities grounded in diffusion models, the United Filter iteratively fine-tunes both state and parameter estimates within a single temporal data assimilation step. Thanks to the robustness of the EnSF, the proposed method offers a promising solution for enhancing our comprehension and modeling of dynamical systems, as demonstrated by results from numerical experiments.

\end{abstract}

\textbf{keyword:} Optimal filtering, state estimation, parameter estimation, diffusion model, ensemble score filter, Bayesian inference


\section{Introduction}

State and parameter estimation present critical challenges in the field of data assimilation. The state estimation problem typically involves a noise-perturbed stochastic dynamical system, where the state may not be directly observable. Instead, partial noisy observations are used to obtain the optimal estimate for the state of the dynamics.
The primary mathematical tool for addressing state estimation is known as  ``optimal filtering''. The goal of the optimal filtering problem is to construct an approximation for the conditional probability density function (PDF) of the target state conditioning on the observational data, referred to as the ``filtering density'', and the optimal state estimation is formulated as the corresponding conditional expectation. Important methods for solving the optimal filtering problem include Bayesian filters such as Kalman-type filters \cite{Kalman1961, UnKF, https://doi.org/10.1029/94JC00572,DataAssimilationUsinganEnsembleKalmanFilterTechnique} and particle filters \cite{particle-filter,MCMC-PF, CT1, Kang-PF,  APF, Sny, vanLeeuwen},, along with other optimal filtering methods related to solving stochastic partial differential equations (SPDEs) \cite{Bao_Zakaid_2015, zakai, Bao_AA20, BaoC20142, Bao_first, BSDE_filter}. Although each of these method has uniquely contributed to the realm of optimal filtering, they all have major drawbacks that cannot be easily addressed.  For example, Kalman filters, replying on the Gaussian assumption for state distributions, can handle relatively high-dimensional problems. But the Gaussian assumption also makes the Kalman filter (and its extended versions) vulnerable when dealing with the nonlinearity that arises in filtering problems. The particle filter, which is also known as a sequential Monte Carlo approach, can alleviate the nonlinearity issue in optimal filtering problems with moderate dimension. However, due to limitations of the Monte Carlo method and challenges associated with high dimensional sampling in Bayesian inference,  extending particle filters to solve high dimensional problems is a formidable task. While the SPDE-related methods are mathematically solid and can provide stable state estimation performance over time, the higher computational cost usually accompanying with SPDE solvers makes the SPDE-related optimal filtering methods less popular in those applications involving extensive computations, such as weather forecasting and state estimation for power systems.

In most studies on state estimation, there is a common assumption that the dynamics of the state are explicitly given. However, practical situations often involve state models with unknown factors. This leads to the challenge of parameter estimation. Methods for parameter estimation can be categorized into the point estimation approach, such like maximum likelihood estimation (MLE)\cite{MLE} and maximum a posteriori estimate (MAP)\cite{Dean2007}, and the Bayesian approach. Approaches to parameter estimation can be broadly categorized into the point estimation approach, including methods like Maximum Likelihood Estimation (MLE) and Maximum A Posteriori Estimate (MAP), and the Bayesian approach. When integrating parameter estimation with state estimation through optimal filtering, a natural choice is to employ a recursive Bayesian inference approach for joint online parameter estimation and dynamic state estimation. The current state-of-the-art method for joint state-parameter estimation is the Augmented Ensemble Kalman Filter (AugEnKF) \cite{Augmentedfilter}. This method considers unknown model parameters as additional state variables and analyzes them alongside other state variables through the ensemble Kalman filter (EnKF). Apparently, the success of the AugEnKF heavily depends on the performance of the EnKF. Moreover, the strategy that uses state estimation to indirectly influence parameter estimation often results in less accurate parameter estimation, as the focus of the augmented approach for joint state-parameter estimation is the state estimation task \cite{Doucet_Stats}.

\vspace{0.5em} 

In this paper, we introduce a \textit{United Filter} method that combines a novel \textit{ensemble score filter (EnSF)} for state estimation with a \textit{direct filter} method for online parameter estimation to solve the joint state-parameter estimation problem. 
The EnSF is a generative method for optimal filtering. It adopts the score-based diffusion model framework to generate an extensive set of samples for the filtering density of the target state variable. As a crucial class of generative machine learning methods, diffusion models utilize noise injection to progressively distort data and then learn to reverse this process for sample generation. Widely applied in image synthesis \cite{NEURIPS2021_49ad23d1,NEURIPS2020_4c5bcfec,NEURIPS2019_3001ef25,DBLP:conf/eccv/CaiYAHBSH20,DBLP:journals/jmlr/HoSCFNS22,DBLP:conf/iclr/MengHSSWZE22}, image denoising \cite{DBLP:conf/iccvw/KawarVE21,DBLP:conf/iccv/LuoH21,NEURIPS2020_4c5bcfec,DBLP:conf/icml/Sohl-DicksteinW15}, image enhancement \cite{DBLP:journals/corr/abs-2112-05149,DBLP:journals/ijon/LiYCCFXLC22,DBLP:journals/pami/SahariaHCSFN23,DBLP:conf/cvpr/WhangDTSDM22} and natural language processing  \cite{DBLP:conf/nips/AustinJHTB21,DBLP:conf/nips/HoogeboomNJFW21,DBLP:journals/corr/abs-2205-14217,DBLP:conf/iclr/SavinovCBEO22,DBLP:conf/icml/YuXMJPGZZW22}, diffusion models rely on the ``score'' to store the information about data. When applying score-based diffusion models to solve the optimal filtering problem, we store the information of filtering densities in score models\cite{Bao2023scorebased}.  In the EnSF, we employ an ensemble approximation for the score model without training it through deep neural networks\cite{Bao2023EnSF}.  With the powerful generative capabilities of diffusion models, it has been demonstrated that EnSF can produce samples from complex target distributions \cite{2023EnSF_Density, song2021scorebased}, and it can track stochastic nonlinear dynamical systems in very high dimensional spaces \cite{Bao2023EnSF}.

The direct filter follows the recursive Bayesian inference framework to solve the online parameter estimation problem. Unlike other recursive Bayesian filter methods, such as the AugEnKF, which augment the state variable with model parameters and use optimal filtering techniques to estimate the augmented state-parameter vector, the direct filter directly utilizes Bayesian inference to ``project'' state dynamics and observational data into the parameter space to identify the optimal estimate for the model parameters \cite{Bao_parameter}. While the direct filter can utilize all the available model and data information to determine the unknown parameters, a notable drawback is its reliance on complete observations of the state variable. In scenarios where the state variable is only partially observed, such as in the joint state-parameter estimation problem, extra state estimation is needed to implement the direct filter.

Leveraging the robust state estimation capabilities of the EnSF, the United Filter method utilizes the estimated state obtained by the EnSF to conduct the direct filter. Through an iterative procedure, the United Filter alternately calibrates state estimation and parameter estimation within a single temporal data assimilation step. 

In the following sections, we provide a brief discussion of the joint state-parameter estimation problem and introduce the recursive Bayesian filter framework in Section \ref{Formulation}. Then, in Section \ref{Derivation}, we derive the United Filter. Numerical experiments will be presented in Section \ref{Numerics}.

\section{Problem setting}\label{Formulation}

\subsection{The joint state-parameter estimation problem}
We consider the following stochastic dynamical system: 
\begin{equation}\label{model:state}
X_{n+1} = f(X_n, \gamma) + \omega_n,
\end{equation}
where $X_{n} \in \mathbb{R}^d$ is the state of the dynamics at a time instant $n$, $f: \mathbb{R}^d \times \mathbb{R}^l \rightarrow \mathbb{R}^d$ is a stochastic dynamical model parameterized by $\gamma \in \mathbb{R}^l$, and $\omega_n \in \mathbb{R}^k$ is a standard $k$-dimensional Gaussian random variable with covariance matrix $\Omega$. 

In many practical scenarios, the true state $X$ within the dynamical system \eqref{model:state} is not directly observable, and the parameter $\gamma$ in the model $f$ is not given either. To estimate the state $X$, as well as the unknown parameter $\gamma$, we acquire partial noisy observations on $X$ through the following observation process
\begin{equation}\label{model:observation}
Y_{n+1} = g(X_{n+1}) + \epsilon_{n+1},
\end{equation}
where $g$ is a nonlinear observation function, and $\epsilon_{n+1}$ is a Gaussian noise with zero mean and a given covariance $\Sigma$.  
The goal of the state-parameter estimation problem for the stochastic systems \eqref{model:state} - \eqref{model:observation} is to find the best estimate for $X_{n+1}$ and $\gamma$ given the observational data $Y_{1:n+1} : = \{Y_i\}_{i=1}^{n+1}$, which contains the information of observations up to time $n+1$.

\subsection{The recursive Bayesian filter}
Since state estimation for $X_{n+1}$ has to be a recursively implemented dynamical procedure, we also formulate parameter estimation problem as a sequential process, and we define the following pseudo parameter process
\begin{equation}\label{Para:dynamics}
\gamma_{n+1} = \gamma_n + \xi_n,
\end{equation}
where $\xi_n$ is a pre-chosen noise term that allows $\gamma$ to explore in the parameter space over time, and $\gamma_{n}$ is our estimate for $\gamma$ at time $n$. Then, we aim to find the ``best'' estimate for $X_{n+1}$ as its conditional expectation, i.e. 
$$\bar{X}_{n+1} = \E[X_{n+1} | Y_{1:n+1}],$$
and we use the estimated state $\bar{X}_{n+1}$ to find an estimate $\bar{\gamma}_{n+1}$ for $\gamma_{n+1}$.

The mathematical technique that computes  $\bar{X}_{n+1}$ is the ``optimal filtering'', and our objective is to find the conditional distribution of the state variable $X_{n+1}$ conditioning on the observations $Y_{1:n+1}$, i.e., $P(X_{n+1} | Y_{1:n+1})$. In the context of the optimal filtering problem, this conditional distribution is often referred to as the ``filtering density''.  It's essential to note that even when the parameter $\gamma$ for the state dynamics in Eq. \eqref{model:state} is known, the state estimation problem remains highly challenging. This difficulty arises from the complexities associated with high-dimensional approximation and Bayesian inference for nonlinear systems. Particularly, in cases where there is a discrepancy between the predicted state distribution (the prior) and the observational data (the likelihood), obtaining the posterior distribution in high-dimensional space becomes a formidable task. Moreover, when the model parameter is unknown, one can expect more substantial discrepancies between the predicted state and the observational data, adding further complexity to the implementation of Bayesian inference.

In this study, we employ the recursive Bayesian filter framework, which is the standard approach for solving the optimal filtering problem, to dynamically estimate the filtering densities for $X_{n+1}$ and $\gamma_{n+1}$. In the subsequent section, we introduce an ensemble score filter method for state estimation and a direct filter method for parameter estimation, and then we provide our United Filter that will effectively combine the state estimation procedure and the parameter estimation procedure.

\vspace{1em}

\section{The United Filter}\label{Derivation}

We first introduce our ensemble score filter in subsection \ref{Section:EnSF}, and then in subsection \ref{Section:DF} we describe a direct filter method to dynamically estimate the unknown parameter in the online manner. In subsection \ref{Section:UF}, we combine the ensemble score filter with the direct filter under the United Filter framework. Our entire algorithm for state-parameter estimation will be summarized in subsection \ref{Summary}.

\subsection{The ensemble score filter for state estimation}\label{Section:EnSF}
In the discussion of the ensemble score filter (EnSF), we begin with the assumption that the parameter $\gamma$ is already known, and therefore our primary focus is on state estimation. Then, after introducing the direct filter method for parameter estimation, we will integrate state estimation with parameter estimation.

In the Bayesian filter, given the filtering density, denoted by $p(X_n | Y_{1:n})$, for $X$ at time instant $n$, we generate the predicted filtering density (the prior), i.e. $p(X_{n+1} | Y_{1:n})$, through the following Chapman-Kolmogorov formula
\begin{equation}\label{Kolmogorov}
p(X_{n+1} | Y_{1:n}) = \int p(X_n | Y_{1:n}) p(X_{n+1} | X_n) dX_n, \hspace{1em} \text{(Prediction)}
\end{equation}
where $p(X_{n+1} | X_n)$ is the transition probability governed by the dynamical model $f$ in Eq. \eqref{model:state}.

Upon reception of the observational data $Y_{n+1}$, we apply the following Bayesian inference procedure to combine the data information (in the form of likelihood) with the prior and get the updated filtering density (the posterior). 
\begin{equation}\label{Bayesian}
p(X_{n+1} | Y_{1:n+1}) \propto p(X_{n+1} | X_n) p(Y_{n+1}|X_{n+1}), \hspace{1em} \text{(Update)}
\end{equation}
where the likelihood function $p(Y_{n+1}|X_{n+1})$ is defined by
\begin{equation}\label{Likelihood}
p(Y_{n+1}|X_{n+1}) \propto \exp\Big[ - \f{1}{2} (g(X_{n+1}) - Y_{n+1})^{\top} \Sigma^{-1} (g(X_{n+1}) - Y_{n+1}) \Big].\hspace{1em} \text{(Likelihood)}
\end{equation}

Through the recursive execution of the prediction procedure as defined in Eq. \eqref{Kolmogorov} and the update procedure as outlined in Eq. \eqref{Bayesian}, one can propagate the filtering density over time. While the prediction procedure can be implemented by running independent samples through the Chapman-Kolmogorov Eq. \eqref{Kolmogorov}, the key challenges in the optimal filtering problem are to accurately characterize the filtering density in high dimensional space and to effectively incorporate new data information into the filtering density. 

\vspace{0.5em}

The main idea of EnSF is to introduce a pair of stochastic differential equations (SDEs), following the score-based diffusion model framework \cite{song2021scorebased}, to establish a connection between a filtering density and the standard normal distribution. This connection ensures that the information of the target filtering density is stored in and represented by the score function. More specifically, we construct two score functions, denoted as $S_{n+1 | n}(\cdot, \cdot)$ and $S_{n+1|n+1}(\cdot, \cdot)$, corresponding to the filtering densities $p(X_{n+1} | Y_{1:n})$ and $p(X_{n+1} | Y_{1:n+1})$, respectively.

In what follows, we briefly recall the score-based diffusion model and explain how the score can connect a target distribution and the standard Gaussian distribution. In a score-based diffusion model, there's a forward stochastic differential equation (SDE) and a reverse-time SDE defined in a pseudo-time domain $t \in [0, T]$, i.e.
\begin{equation}\label{DM:forward}
d Z_{t} = b(t) Z_{t} dt + \sigma(t) dW_{t}, \qquad \text{(Forward SDE)}
\end{equation}
\begin{equation}\label{DM:reverse}
d Z_{t} = \big[ b(t) Z_{t} - \sigma^2(t) S(Z_{t}, t)\big]dt + \sigma(t) d\ola{W}_{t}, \qquad \text{(Reverse-time SDE)}
\end{equation}
where $W_{t}$ is a standard $d$-dimensional Brownian motion, $\int \cdot dW_{t}$ denotes a standard It\^{o} type stochastic integral, $\int \cdot d\ola{W}_{t}$ is a backward It\^{o} integral, and $b$ and $\sigma$ in Eq. \eqref{DM:forward} are the drift coefficient and the diffusion coefficient, respectively. In this work, we let the terminal time $T$ in the diffusion model be $T=1$. The function $S(\cdot, \cdot)$ in Eq. \eqref{DM:reverse}, defined by
$$S(z, t) : = \nabla \log(Q_{t}(z)), \quad z \in \mathbb{R}^d,$$
is the so-called ``score'', where $Q_{t}$ is the probability distribution of solution $Z_{t}$ of the forward SDE \eqref{DM:forward}, which transforms any target distribution $Q_0$ of $Z_0$ to the standard Gaussian distribution $N(0, I_d)$ with properly chosen drift $b$ and diffusion $\sigma$ in the diffusion model system.

Given explicit definitions of $b$ and $\sigma$, the information of the initial distribution of $Z_0$ is stored in score $S$. In fact, the score function establishes a connection with the distribution $Q_0(Z_0)$ in the following manner:
\begin{equation}\label{Score:analytic}
\begin{aligned}
S(Z_{t}, t) = & \nabla \log (Q_{t}(Z_{t}))= \nabla \log \Big( \int_{\mathbb{R}^d} Q_{t}(Z_{t} | Z_0) Q_0(Z_0) dZ_0 \Big) \\
= & \int_{\mathbb{R}^d} - \f{Z_{t} - \alpha_{t} Z_0}{\beta^2_{t}} w_{t}(Z_{t}, Z_0) Q_0(Z_0) dZ_0,
\end{aligned}
\end{equation} 
where 
$$w_{t}(Z_{t}, Z_0):= \f{Q_{t}(Z_{t} | Z_0)}{\int_{\mathbb{R}^d} Q_{t}(Z_{t} | Z'_0) Q_{0}(Z'_0) dZ'_0}$$
is a weight factor. 

In this work, we let
$$b(t) = \f{d\log \alpha_{t}}{dt}, \qquad \sigma^2(t) = \f{d \beta^2_{t}}{d t} - 2 \f{d \log \alpha_{t}}{d t} \beta^2_{t},$$
with
$\alpha_{t} = 1 - t$ and $\beta_{t} = t$ for $t \in [0, 1]$. Therefore, for any given initial random variable $Z_0$ (with the target distribution  $Q_0$), the conditional PDF $Q_{t}(Z_{t} | Z_0)$ can be described by the following expression
$$Q_{t}(Z_{t} | Z_0) = N(\alpha_{t} Z_0, \beta^2_{t} I_{d}).$$
Clearly, when $t \rightarrow 1$, the conditional distribution $Q_{t}(Z_{t} | Z_0)$ approaches $N(0, I_d)$, and as $t \rightarrow 0$, $Q_{t}(Z_{t} | Z_0)$ approaches the distribution of $Z_0$.

Based on the above discussions, we define scores $S_{n+1 | n}$ and $S_{n+1 | n+1}$ corresponding to the prior filtering density $p(X_{n+1} | Y_{1:n})$ and the posterior filtering density $p(X_{n+1} | Y_{1:n+1})$, respectively, as follows:
$$\text{Prior filtering score: } S_{n+1 | n}: \text{set $Z_0 = X_{n+1}|Y_{1:n}$ and hence $Q_{0}(Z_0) = p(X_{n+1} | Y_{1:n})$ },$$
$$\text{Posterior filtering score: } S_{n+1 | n+1}: \text{set $Z_0 = X_{n+1}|Y_{1:n+1}$ and hence $Q_{0}(Z_0) = p(X_{n+1} | Y_{1:n+1})$}.$$ 

In order to obtain $S_{n+1 | n}$, data samples for $Z_0 = X_{n+1}|Y_{1:n}$ are needed. To this send, we first generate samples $\{x_{n|n}^j\}_{j=1}^J$ from the posterior filtering score $\bar{S}_{n|n}$ at time instant $n$, where $\bar{S}_{n|n}$ is an approximation for $S_{n|n}$. Then, we simulate samples $\{x_{n|n}^j\}_{j=1}^J$ through the state model \eqref{model:state} $f$ to get $\{x_{n+1|n}^j\}_{j=1}^J$, i.e.,
\begin{equation}\label{samples:predict}
x_{n+1|n}^j = f(x_{n|n}^j, \gamma, \omega_{n, j}), \quad j = 1, 2, \cdots, J,
\end{equation}
and $\{x_{n+1|n}^j\}_{j=1}^J$ are the desired samples for $Z_0 = X_{n+1}|Y_{1:n}$.

Then, the remaining challenge lies in approximating the score $S_{n+1 | n}$. The standard approach to approximate score functions involves the use of  deep neural networks \cite{Bao2023scorebased}. In a recent study, we introduced a training-free ensemble approach for approximating the score function, and this approach will be employed in this work.

Following the analytic expression for the score model that we provided in Eq. \eqref{Score:analytic}, we approximate $S_{n+1 | n}$ at a given $z\in\mathbb{R}^d$ and $t \in [0, 1]$ as follows:
\begin{equation}\label{Approx:S-prior}
S_{n+1|n}(z, t) \approx \bar{S}_{n+1|n}(z, t) := \sum_{m=1}^{M} - \f{z - \alpha_{t} \cdot x_{n+1|n}^{j_m}}{\beta^2_{t}} \bar{w}_{t}( z, x_{n+1|n}^{j_m} ),
\end{equation}
where $\{x_{n+1|n}^{j_m}\}_{m=1}^M$ is a mini-batch of samples in the sample set  $\{x_{n+1|n}^j\}_{j=1}^J$ introduced by Eq. \eqref{samples:predict}, and $\bar{w}_{t}$ is an approximation to weight $w_{t}$, which is defined as follows:
$$\bar{w}_{t}( z, x_{n+1|n}^{j'_m}) : = \f{Q_{t}(z | x_{n+1|n}^{j'_m} )}{\sum_{m=1}^M Q_{t}(z | x_{n+1|n}^{j_m} )}.$$

In the process of approximating the score $S_{n+1|n+1}$, which is associated with the posterior filtering density $p(X_{n+1} | Y_{1:n+1})$, we face a challenge due to the absence of samples from the posterior filtering density.  To address this, we analytically incorporate the likelihood information into the current estimated score $\bar{S}_{n+1|n}$, which is corresponding to the prior filtering density $p(X_{n+1} | Y_{1:n})$. 

Specifically,  we approximate $S_{n+1 | n+1}$ as follows:
\begin{equation}\label{Approx:S-posterior}
\bar{S}_{n+1 | n+1}(z, t) := \bar{S}_{n+1|n}(z, t) + h(t) \nabla \log p(Y_{n+1} |  z),
\end{equation}
where $\bar{S}_{n+1|n}$ is the training-free approximation introduced in Eq. \eqref{Approx:S-prior}, $p(Y_{n+1} |  \cdot )$ is the likelihood given by Eq. \eqref{Likelihood}, and $h$ is a damping function satisfying
$$\text{$h(\cdot)$ monotonically decreases in $[0, 1]$ with $h(0) = 1$ and $h(1) = 0$}.$$
Note that, by taking the gradient to the log of the posterior filtering density, we have
\begin{equation}
\nabla \log p(X_{n+1} | Y_{1:n+1}) = \nabla \log p(X_{n+1} | Y_{1:n}) + \nabla \log p(Y_{n+1} | X_{n+1}).
\end{equation}
Moreover, from the definition of the score function, we know that 
$$S_{n+1|n}(Z_0, 0) = \nabla \log p(X_{n+1} | Y_{1:n}).$$
Therefore, when $t = 0$ in Eq. \eqref{Approx:S-posterior}, 
$$\bar{S}_{n+1 | n+1}(Z_0, 0) =  \nabla \log p(X_{n+1} | Y_{1:n}) + \nabla \log p(Y_{n+1} |  X_{n+1}) = \nabla \log p(X_{n+1} | Y_{1:n+1}),$$
which targets on the recovery of the posterior filtering density. 

On the other hand, at pseudo-time $t = 1$, $\bar{S}_{n+1|n+1}(z, 1) = \bar{S}_{n+1|n}(z, 1)$ as given in Eq.~\eqref{Approx:S-posterior}, which is corresponding to the standard Gaussian distribution $N(0, I_d)$. In this way, our approximated score $\bar{S}_{n+1 | n+1}$ does connect the standard Gaussian distribution to the posterior filtering density $p(X_{n+1} | Y_{1:n+1})$ in the reverse-time diffusion procedure. We refer to \cite{Bao2023EnSF} for more detailed discussions on the approximation of posterior score $S_{n+1 | n+1}$. In this work, we let $h(t) = 1 - t$ in our numerical experiments in Section \ref{Numerics}.

With approximated score $\bar{S}_{n+1|n+1}$ in Eq. \eqref{Approx:S-posterior}, we can generate a set of samples $\{x_{n+1|n+1}^j\}_{j=1}^J$ that follow the posterior filtering density $p(X_{n+1} | Y_{1:n+1})$ from the Gaussian distribution $N(0, I_d)$ through the reverse-time SDE \eqref{DM:reverse}, and the estimated state can be obtained by 
$$\bar{X}_{n+1} \approx \f{1}{J}\sum_{j=1}^{J} x_{n+1|n+1}^j.$$

Note that the above EnSF framework requires the knowledge of the model parameter $\gamma$. In what follows, we shall provide a direct filter method to dynamically estimate the unknown parameter in model $f$ upon reception of observational data.

\subsection{A direct filter method for online parameter estimation}\label{Section:DF}

To introduce the direct filter, we assume that the estimates $\{\bar{X}_{n}\}_{n}$ for the state variable $\{X_{n}\}_{n}$ are available, and we shall apply the estimated state to derive an estimate for the unknown parameter. 

The primary concept behind the direct filter is to ``project'' the information about the state variable onto the parameter space through Bayesian inference and determine the best estimate for the unknown parameter. With the estimated states $\{\bar{X}_{n}\}_{n}$, we let the conditional expectation 
\begin{equation*}
\bar{\gamma}_{n+1} = \E[\gamma_{n+1} | \bar{X}_{1:n+1}]
\end{equation*}
be our estimate for $\gamma_{n+1}$. Note that estimates $\{\bar{X}_{n}\}_{n}$ for $\{X_{n}\}_{n}$ are conditioning on the observations $\{Y_{n}\}_n$, and the optimal filtering procedure for finding $\{\bar{X}_{n}\}_{n}$ also incorporates the dynamical model $f$. Therefore, the conditional expectation $\E[\gamma_{n+1} | \bar{X}_{1:n+1}]$, as an estimate for $\gamma_{n+1}$, has considered both the state dynamics and the observational data.

To achieve this, we target on generating the parameter distribution $p(\gamma_{n+1}|\bar{X}_{1:n+1})$ for $\gamma_{n+1}$.
More specifically, the prior parameter distribution $p(\gamma_{n+1}|\bar{X}_{1:n})$ is produced through the pseudo parameter dynamics \eqref{Para:dynamics} and the previous estimated parameter distribution $p(\gamma_{n}|\bar{X}_{1:n})$. Then, we adopt the following Bayesian formula to obtain our description for the estimated parameter at time instant $n+1$
\begin{equation}\label{Bayes:para}
p(\gamma_{n+1}|\bar{X}_{1:n+1})\propto p(\gamma_{n+1}|\bar{X}_{1:n}) p(\bar{X}_{n+1}|\gamma_{n+1}),
\end{equation}
where $p(\bar{X}_{n+1}|\gamma_{n+1})$ is the likelihood that measures the discrepancy between the predicted parameter and our estimated state $\bar{X}_{n+1}$.

In this work, we utilize a particle implementation \cite{Bao_Cogan20, Bao_parameter} to execute the parameter Bayesian \eqref{Bayes:para}\footnote{Other recursive Bayesian inference techniques can also be applied to implement Eq. \eqref{Bayes:para}.}. We assume that we have a set of parameter particles, denoted by $\{\gamma^{k}_n\}_{k=1}^K$, from the parameter estimation distribution  $p(\gamma_{n}|\bar{X}_{1:n})$, i.e.,
$$p(\gamma_{n}|\bar{X}_{1:n}) \approx \f{1}{K}\sum_{k=1}^K \delta_{\gamma^{k}_n}(\gamma_{n}),$$
where $K$ is a pre-determined integer representing the number of particles that we choose to approximate the parameter distribution, and $\delta_{\cdot}$ is the Dirac delta function. Then, we add random perturbations to the parameter particle cloud (as discussed in \eqref{Para:dynamics}) to get a set of predicted particles $\{\tilde{\gamma}^{k}_{n+1}\}_{k=1}^K$, i.e.
\begin{equation}
\tilde{\gamma}^{k}_{n+1} = \gamma^{k}_n + \xi_n^k,
\end{equation}
where $\xi_n^k \sim N(0, \Gamma)$ is a Gaussian type random variable with a pre-determined covariance matrix $\Gamma$ representing the exploration range for the unknown parameter.

With a set of predicted parameter particles $\{\tilde{\gamma}^{k}_{n+1}\}_{k=1}^K$, we evaluate the likelihood of each particle by comparing it with the estimated state $\bar{X}_{n+1}$ through the likelihood $p(\bar{X}_{n+1}|\gamma_{n+1})$. Here $\bar{X}_{n+1}$ is obtained by the EnSF for state estimation, and it contains information from the observational data $Y_{n+1}$ that reflects the actual model parameter. 

However, the model parameters are not directly observable. In the direct filter method, we input each parameter particle $\gamma_{n+1}^k$ into the state model $f$ to obtain a parameter-dependent state $\bar{X}^{\gamma_{n+1}^k}_{n+1} = f(\bar{X}_{n}, \gamma_{n+1}^k) + \omega_n^{k}$, and we use the following likelihood function
\begin{equation}\label{gamma:likelihood}
p(\bar{X}_{n+1}|\tilde{\gamma}^k_{n+1}) \propto \exp\Big[ - \f{1}{2} \big(\bar{X}^{\tilde{\gamma}_{n+1}^k}_{n+1} - \bar{X}_{n+1}\big)^{\top} \Omega^{-1} \big(\bar{X}^{\tilde{\gamma}_{n+1}^k}_{n+1} - \bar{X}_{n+1}\big) \Big]
\end{equation}
to compare $\bar{X}^{\gamma_{n+1}^k}_{n+1}$ with the estimated state $\bar{X}_{n+1}$. 

In this way, the likelihood $p(\bar{X}_{n+1}|\tilde{\gamma}^k_{n+1})$ assigns a weight to the parameter particle $\tilde{\gamma}^k_{n+1}$, and we describe the conditional parameter distribution, i.e., $p(\gamma_{n+1}|\bar{X}_{1:n+1})$ introduced in \eqref{Bayes:para}, by using the following weighted sampling scheme
$$p(\gamma_{n+1}|\bar{X}_{n+1}) \approx  \sum_{k=1}^K w_k \delta_{\tilde{\gamma}^{k}_{n+1}}(\gamma_{n+1}),$$
where the weight $w_k \propto p(\bar{X}_{n+1}|\tilde{\gamma}^k_{n+1})$. 

To mitigate the degeneracy issue often encountered in particle filter approaches, we implement a resampling procedure. This process involves creating additional duplicates of particles in $\{\tilde{\gamma}^{k}_{n+1}\}_{k=1}^K$ with higher weights, resulting in a set of equally weighted particles denoted as $\{\gamma^{k}_{n+1}\}_{k=1}^K$\cite{particle-filter}. Then, our estimate for $\gamma_{n+1}$ becomes
\begin{equation}\label{DF:estimate}
\bar{\gamma}_{n+1} = \f{1}{K}\sum_{k=1}^K \gamma^k_{n+1}.
\end{equation}

\subsection{An iterative procedure that unites the ensemble score filter and the direct filter}\label{Section:UF}
With the EnSF method for state estimation (introduced in Section \ref{Section:EnSF}) and the direct filter method for parameter estimation (introduced in Section \ref{Section:DF}), in this subsection we introduce an iterative procedure to combine the EnSF and the direct filter.

At time instant $n$, we assume that we have an estimated state $\bar{X}_n$, accompanied by a posterior score $\bar{S}_{n|n}$ corresponding to the posterior filtering density $p(X_n|Y_{1:n})$. Additionally, we assume that a set of parameter particles $\{\gamma^{k}_{n}\}_{k=1}^K$ is available to describe the parameter $\gamma_n$, and we denote the mean value of $\{\gamma^{k}_{n}\}_{k=1}^K$ as $\bar{\gamma}_n$. The United Filter algorithm for joint state-parameter estimation at time instant $n+1$ comprises two iterative stages: (I) state estimation with the current estimated parameter, and (II) parameter estimation with the current estimated state. The iteration procedure, indexed by $l=0, 1, 2, \cdots, L$, begins with the initial estimate for $\gamma_{n+1}$ set as $\bar{\gamma}^{(0)}_{n+1} = \bar{\gamma}_n$ with parameter particles chosen as $\{\gamma^{k, (0)}_{n+1} \}_{k=1}^K : = \{\gamma^{k}_{n} \}_{k=1}^K$.

\subsubsection*{\textbf{(I). State estimation to obtain $\bar{X}^{(l)}_{n+1}$ with the current estimated parameter $\bar{\gamma}^{(l)}_{n+1}$}}
At iteration stage $l$ for the time period $n$ to $n+1$, we first generate a set of $J$ state samples $\{x_{n|n}^j\}_{j=1}^J$ from the reverse-time SDE \eqref{DM:reverse} with score model $\bar{S}_{n|n}$. Then, we propagate $\{x_{n|n}^j\}_{j=1}^J$ through the state dynamics \eqref{model:state} based on the current estimated parameter $\bar{\gamma}^{(l)}_{n+1}$ to obtain a set of predicted state samples $\{x^{(l), j}_{n+1|n}\}_{j=1}^J$ for the prior filtering density, i.e.,
\begin{equation}\label{parameter2sample}
x^{(l), j}_{n+1|n} = f(x^j_{n|n}, \bar{\gamma}^{(l)}_{n+1}) + \omega^{(l)}_{n, j}.
\end{equation}
Next, we compute the estimated prior score $\bar{S}^{(l)}_{n+1|n}$ using the training-free scheme \eqref{Approx:S-prior}. This computation involves using data samples $\{x^{(l), j}_{n+1|n}\}_{j=1}^J$, where the current estimated model parameter, i.e., $\gamma^{(l)}_{n+1}$, is already incorporated into $\{x^{(l), j}_{n+1|n}\}_{j=1}^J$ through Eq. \eqref{parameter2sample}. 

To integrate the observational data $Y_{n+1}$ in state estimation, we approximate the posterior score $\bar{S}^{(l)}_{n+1|n+1}$ by using scheme \eqref{Approx:S-posterior}. As a result, we can generate (posterior) state samples $\{x^{(l), j}_{n+1|n+1}\}_{j=1}^J$ from score $\bar{S}^{(l)}_{n+1|n+1}$ to evaluate $X_{n+1}$, and we denote our updated estimate as $\bar{X}^{(l)}_{n+1}$.

\subsubsection*{\textbf{(II). Parameter estimation to obtain $\bar{\gamma}^{(l+1)}_{n+1}$ with the current estimated state $\bar{X}^{(l)}_{n+1}$ }}
The direct filter based parameter estimation relies on the current estimate $\bar{X}^{(l)}_{n+1}$ for the state variable $X_{n+1}$, and we slightly modify the direct filter method introduced in Section \ref{Section:DF} by replacing $\bar{X}_{n+1}$ with $\bar{X}^{(l)}_{n+1}$. In the practical numerical implementation, at the $l+1$-th iteration stage, we generate a set of predicted parameter particles as 
$$\tilde{\gamma}_{n+1}^{k, (l+1)} = \gamma_{n+1}^{k, (l)} + + \xi_n^{k, (l+1)}, \quad k = 1, 2, \cdots, K,$$
and we re-write the likelihood in Eq. \eqref{gamma:likelihood} as 
\begin{equation}\label{gamma:likelihood:iteration}
p(\bar{X}^{(l)}_{n+1}|\tilde{\gamma}^{k, (l+1)}_{n+1}) \propto \exp\Big[ - \f{1}{2} \big(\bar{X}^{\tilde{\gamma}_{n+1}^{k, (l+1)}}_{n+1} - \bar{X}^{(l)}_{n+1}\big)^{\top} \Omega^{-1} \big(\bar{X}^{\tilde{\gamma}_{n+1}^{k, (l+1)}}_{n+1} - \bar{X}^{(l)}_{n+1}\big) \Big]
\end{equation}
to compare the parameter particle based state prediction $\bar{X}^{\tilde{\gamma}_{n+1}^{k, (l+1)}}_{n+1}$ with the EnSF estimated state $\bar{X}^{(l)}_{n+1}$, and the likelihood $p(\bar{X}^{(l)}_{n+1}|\tilde{\gamma}^{k, (l+1)}_{n+1})$ also serves as the weight for the parameter particle $\tilde{\gamma}^{k, (l+1)}_{n+1}$. 

After a resampling procedure, we obtain a set of equally-weighted parameter particles $\{\gamma^{k, (l+1)}_{n+1} \}_{k=1}^K$, and the direct filter estimated parameter at the $l+1$-th iteration stage is $\bar{\gamma}^{(l+1)}_{n+1} = \f{1}{K}\sum_{k=1}^K \gamma^{k, (l+1)}_{n+1}$.

\vspace{2em}

We iteratively perform the above state estimation stage (I) and parameter estimation stage (II). Following $L$ iterations, we consider that the observational data $Y_{n+1}$ to be sufficiently assimilated into both the estimated state and the estimated parameter. Subsequently, we define the estimated state at time instant $n+1$ as 
$$\bar{X}_{n+1} := \bar{X}^{(L)}_{n+1},$$ 
and the estimated parameter as 
$$\bar{\gamma}_{n+1} := \bar{\gamma}^{(L)}_{n+1}.$$

\subsection{Summary of the algorithm for the United Filter}\label{Summary}
We summarize the entire United Filter algorithm in Algorithm 1.

\noindent\makebox[\linewidth]{\rule{\textwidth}{0.5pt}}\\
\vspace{-0.5cm}
\newline {\bf Algorithm 1: the pseudo-algorithm for the United Filter}\vspace{-0.2cm} \\ 
\noindent\makebox[\linewidth]{\rule{\textwidth}{0.5pt}}
\vspace{-0.3cm}
\newline1:\, {\bf Input}:  the state equation $f(X, \gamma, \omega)$, the observation function $g(X)$, the prior density $P(X_0)$, and the initial parameter particles $\{\gamma^{k}_0\}_{k=1}^K$;
\vspace{0.1cm}
\newline2: Generate $J$ samples $\{x_{0|0}^j\}_{j=1}^J$ from the prior $P(X_0)$;
\vspace{0.1cm}
\newline3: \,{\bf for} $n = 0, \ldots, $
\vspace{0.1cm}
\newline4: \qquad Set $\bar{\gamma}^{(0)}_{n+1} = \bar{\gamma}_{n}$ and $\{\gamma^{k, (0)}_{n+1}\}_{k=1}^{K} = \{\gamma^{k}_{n}\}_{k=1}^{K} $;
\vspace{0.1cm}
\newline5: \qquad {\bf for} $l = 0, 1, 2, \cdots, L-1 $
\vspace{0.1cm}
\newline6: \qquad \quad\; - Carry out the EnSF based state estimation by using the estimated parameter $\bar{\gamma}^{(l)}_{n+1}$ and obtain the estimated state $\bar{X}^{(l)}_{n+1}$ as well as $\bar{S}^{(l)}_{n+1|n+1}$;
\newline7: \qquad \quad\; - Carry out  the direct filter based parameter estimation by using the estimate state $\bar{X}^{(l)}_{n+1}$ and obtain the estimated parameter $\bar{\gamma}^{(l+1)}_{n+1}$ as well as parameter particles $\{\gamma^{k, (l+1)}_{n+1}\}_{k=1}^{K} $;
\newline11: \quad\; {\bf end}\vspace{0.1cm} 
\newline10: \quad\;  - Carry out the EnSF based state estimation by using the estimated parameter $\bar{\gamma}^{L}_{n+1}$ to obtain $\bar{X}^{(L)}_{n+1}$ as well as $\bar{S}^{(L)}_{n+1|n+1}$; 
\newline10: \quad\; Set $\bar{\gamma}_{n+1} = \bar{\gamma}^{(L)}_{n+1}$, $\{\gamma^{k}_{n+1}\}_{k=1}^{K} = \{\gamma^{k, (L)}_{n+1}\}_{k=1}^{K}$, $\bar{X}_{n+1} = \bar{X}^{(L)}_{n+1}$, and $\bar{S}_{n+1|n+1} = \bar{S}^{(L)}_{n+1|n+1}$;
\vspace{0.1cm}
\newline11: {\bf end}\vspace{-0.1cm} \\
\noindent\makebox[\linewidth]{\rule{\textwidth}{0.5pt}}

\vspace{2em}

\section{Numerical experiments}\label{Numerics}
We present three numerical examples to showcase the effectiveness of our United Filter in solving joint state-parameter estimation problems. In the first example, we examine a noise-perturbed SIR model, commonly known as the epidemic model, to demonstrate the baseline performance of the United Filter in estimating both state and model parameters. For the second example, we tackle a high-dimensional state-parameter estimation problem related to a finite difference scheme-discretized Fokker-Planck equation -- an essential model in partial differential equations (PDEs). In this case, our goal is to utilize partial observational data for estimating the solution of the Fokker-Planck equation as a high-dimensional state vector, considering unknown model parameters influencing the equation's propagation. In the last example, we consider a high-dimensional Lorenz-96 model with three unknown model parameters within the system. Simultaneously, we track the state of a Lorenz-96 system in the $200$-dimensional space. The chaotic nature of Lorenz systems adds substantial complexities to this joint state-parameter estimation problem. Moreover, inaccurate estimates for model parameters can lead to significant deviations in state estimation, posing even more formidable challenges in high-dimensional approximation. In the third example, we shall demonstrate the superior performance of our United Filter by comparing it with the state-of-the-art approach for the joint state-parameter estimation problem, i.e., the Augmented Ensemble Kalman Filter.

\subsection{Example 1: An epidemic model}

We first consider the SIR epidemic model:
\begin{equation}\label{SIR:model}
\begin{aligned}
S_{n+1} =& S_{n}  - B S_n I_n \Delta t + \sigma \xi^{1}_n, \\
I_{n+1} =& I_n + \big( B S_n I_n - K I_n \big) \Delta t + \sigma \xi^{2}_n, \\
R_{n+1} =& R_n + K I_n \Delta t + \sigma \xi^{3}_n,
\end{aligned}
\end{equation}
where $S$ is the susceptible fraction of the population, $I$ is the infected fraction of the population, $R$ is the recovered fraction of the population, and $\xi^{i}_n \sim N(0, I)$ is a standard Gaussian noise that perturbs the SIR model with coefficient $\sigma$. There are two parameters in mode \eqref{SIR:model}: $B$ and $K$. The parameter $B$ stands for the contact rate, which reflects the average number of contacts per person in the time unit $\Delta t$, and $K$ is the recovery rate, which is the probability of an infectious individual recovering in the time interval $\Delta t$. The goal of the state-parameter estimation problem for the epidemic model is to estimate the state  $X = [S, I, R]^{\top}$, as well as the unknown contact rate $B$ and recovery rate $K$. 

We assume that one can obtain direct observations on $S$, $I$, $R$, and we let the observation be
$$Y_{n+1} = [S_{n+1}, I_{n+1}, R_{n+1}]^{\top} + \delta \zeta_{n+1}, $$
where $\zeta_{n+1} \sim (0, I_3)$ is a 3-dimensional standard Gaussian random variable representing the observational noise, and $\delta$ is the size of noise that reflects the accuracy of our measurements on $[S_{n+1}, I_{n+1}, R_{n+1}]^{\top}$.  

\begin{figure}[h!]
\begin{center}
\includegraphics[scale = 0.25]{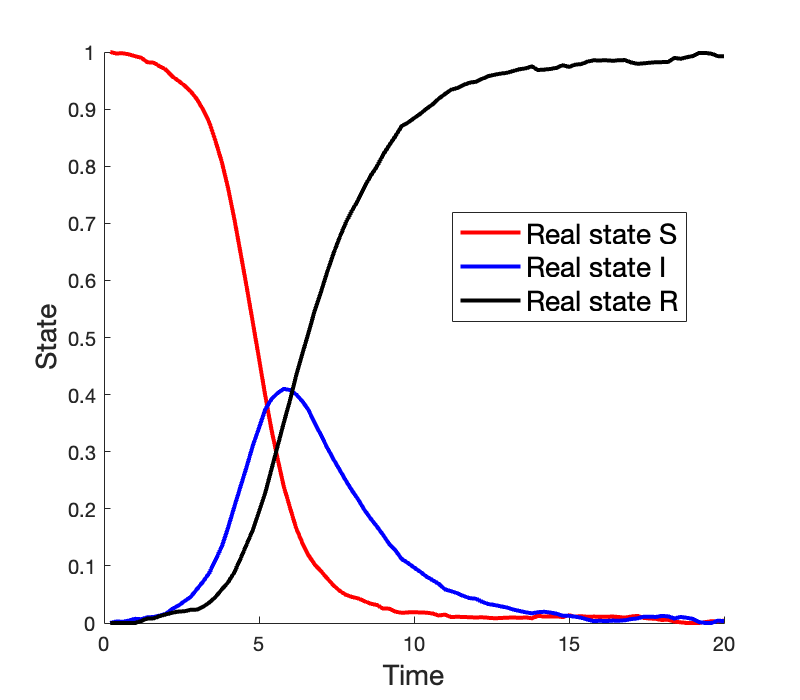} 
\end{center} \vspace{-0.5em}
\caption{Example 1. Real state evolution of the SIR model.}\label{Ex1:Epidemic_real}
\end{figure}

In our numerical experiments, we let the true parameter values be $B = 0.5$ and $K = 2$, and we let  $\sigma = 0.005 \sqrt{\Delta t}$. The initial state is chosen as $S_0 = 1$, $I_0 = 10^{-6}$, and $R_0 = 0$.  In Figure \ref{Ex1:Epidemic_real}, we present the true state of $S$, $I$, and $R$ over the time interval $[0, 20]$, which is discretized into $100$ steps with time step-size $\Delta t = 0.2$.  

\begin{figure}[h!]
\begin{center}
\subfloat[Parameter estimation: B]{\includegraphics[scale = 0.23]{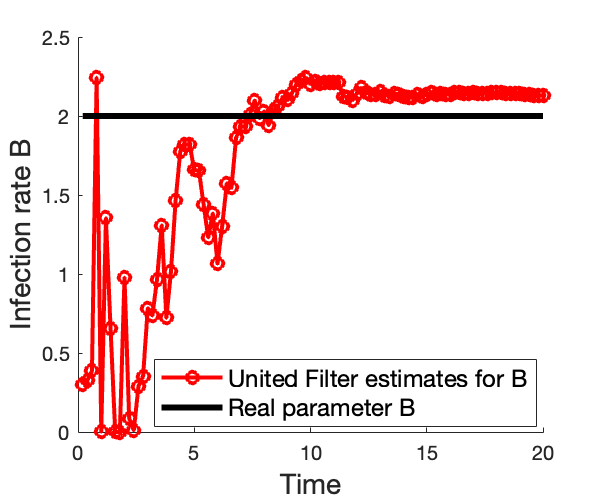} }  
\subfloat[Parameter estimation: K]{\includegraphics[scale = 0.23]{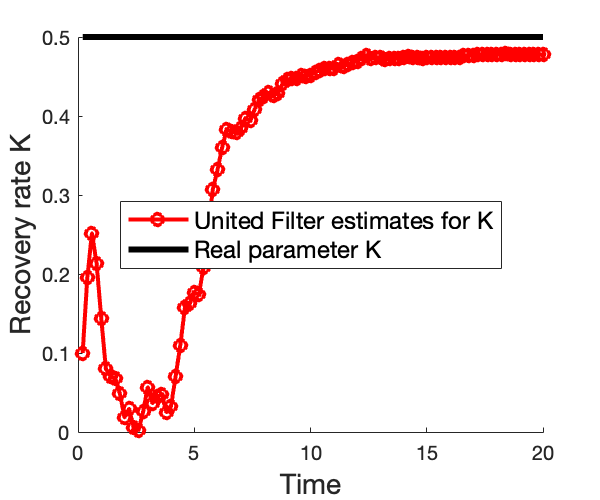} }  
\subfloat[State estimation]{\includegraphics[scale = 0.18]{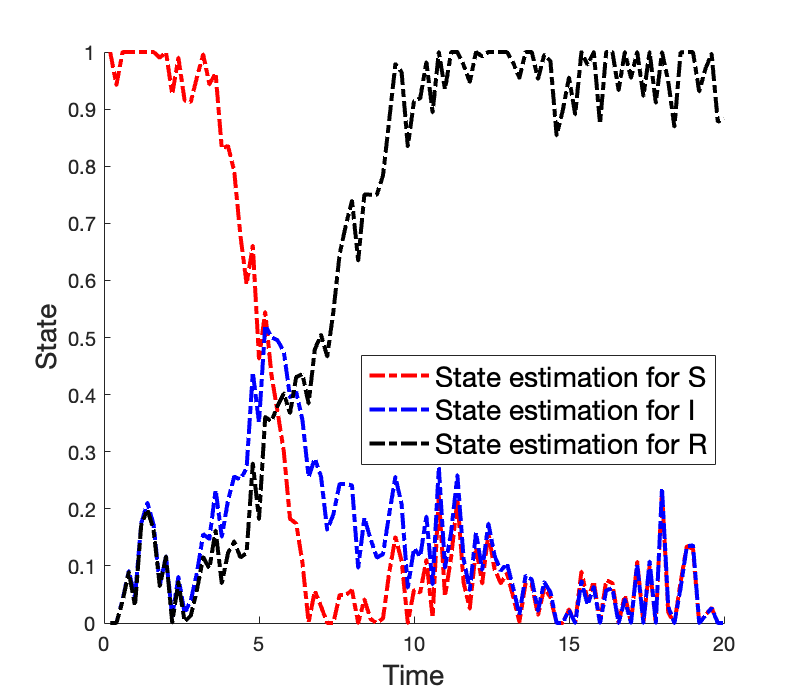} } 
\end{center} \vspace{-0.5em}
\caption{Example 1. State - parameter estimation with observation noise level  $\delta = 0.1$. }\label{Ex1:Epidemic_01}
\end{figure}

A practical issue that we consider in this state-parameter estimation problem for the epidemic model is the accuracy of measurements. In this example, we incrementally increase the noise added to the observations from $\delta = 0.01$ to $\delta = 0.05$, and then to $\delta = 0.1$. Note that when $\delta = 0.1$, our measurements may potentially encompass errors of $10\%$ in the total population. 
\begin{figure}[h!]
\begin{center}
\subfloat[Parameter estimation: B]{\includegraphics[scale = 0.23]{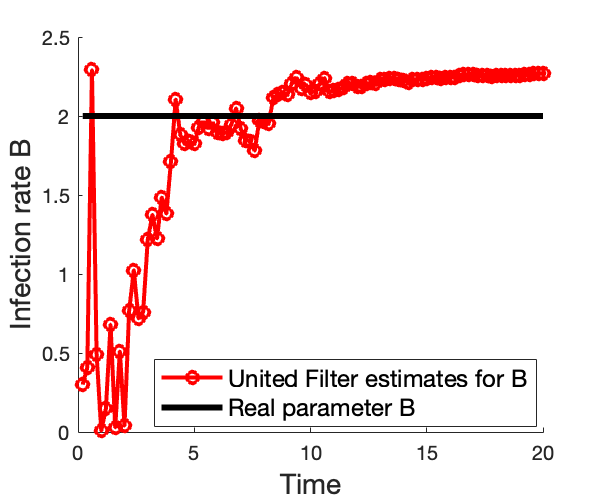} }   
\subfloat[Parameter estimation: K]{\includegraphics[scale = 0.23]{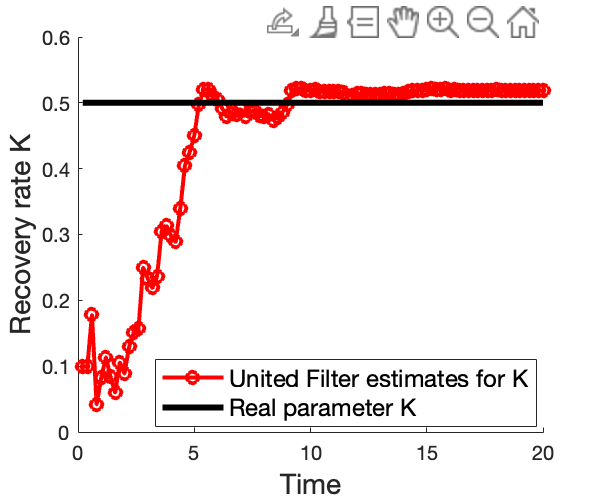} }  
\subfloat[State estimation]{\includegraphics[scale = 0.18]{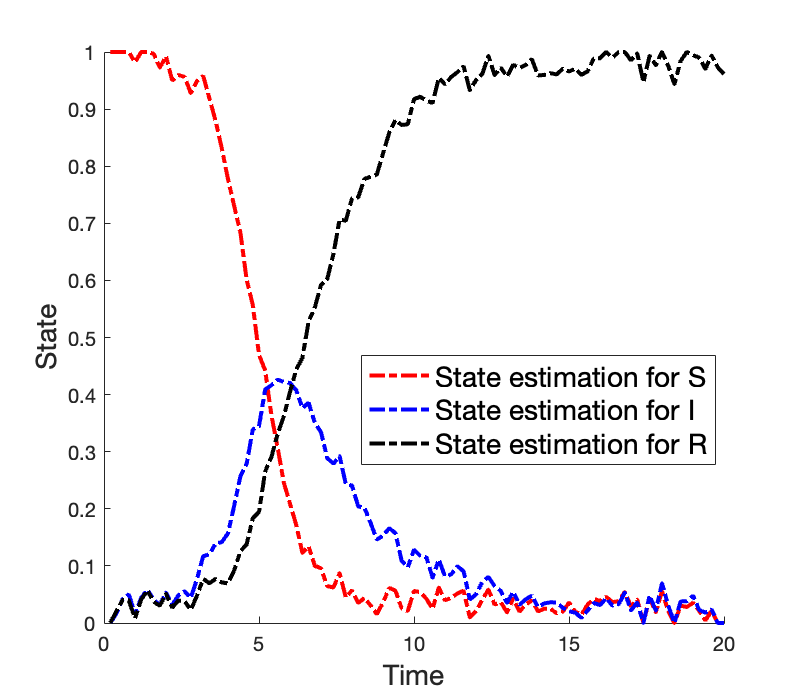} } 
\end{center} \vspace{-0.5em}
\caption{Example 1. State - parameter estimation with observation noise level  $\delta = 0.01$. }\label{Ex1:Epidemic_001}
\end{figure}
In Figure \ref{Ex1:Epidemic_001}, we let the level of observation noise be $\delta = 0.01$ and present the performance of the United Filter in estimating parameter $B$ (the contact rate), parameter $K$ (recovery rate), and the state $S$ (susceptible fraction), $I$ (infected fraction), and $R$ (recovered fraction). 
\begin{figure}[h!]
\begin{center}
\subfloat[Parameter estimation: B]{\includegraphics[scale = 0.23]{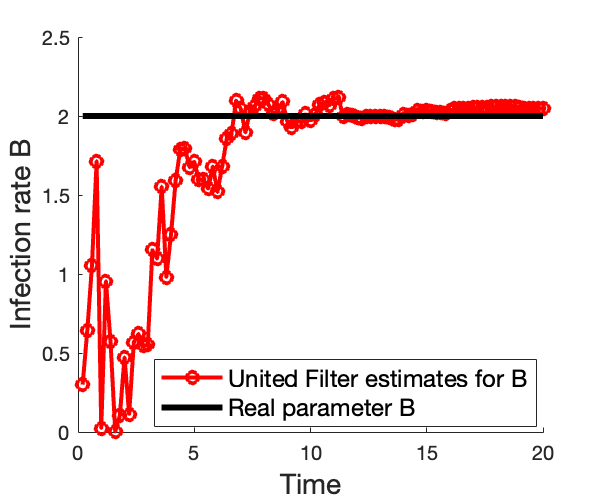} }  
\subfloat[Parameter estimation: K]{\includegraphics[scale = 0.23]{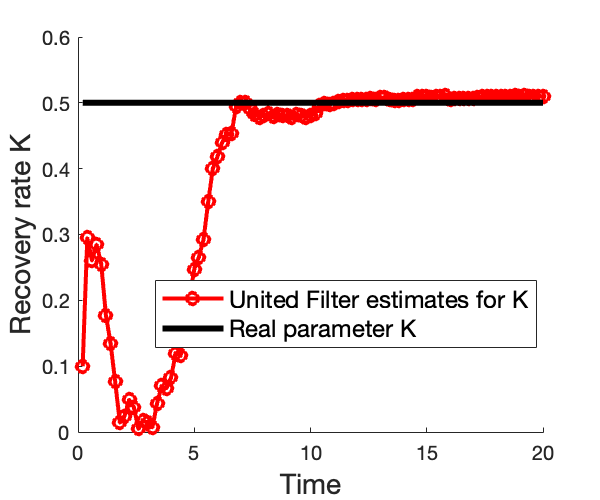} }  
\subfloat[State estimation]{\includegraphics[scale = 0.18]{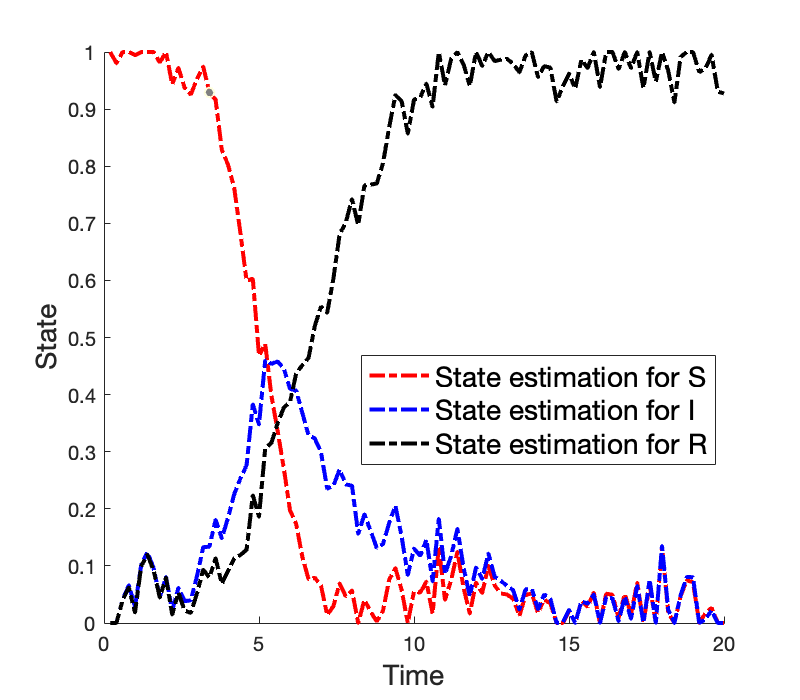} } 
\end{center} \vspace{-0.5em}
\caption{Example 1. State - parameter estimation with observation noise level  $\delta = 0.05$. }\label{Ex1:Epidemic_005}
\end{figure}
Similarly, in Figure \ref{Ex1:Epidemic_005} and Figure \ref{Ex1:Epidemic_01}, we demonstrate the performance of the United Filter in parameter estimation and state estimation by setting the observation noise level at $\delta = 0.05$ and $\delta = 0.1$, respectively. From these figures, we can see that our parameter estimation results remain reliable, regardless of the varying levels of observational noise. However, the accuracy of state estimation is closely tied to the quality of observations.

\subsection{Example 2: Fokker-Planck equation}
Next, we solve a more challenging problem. In the second example, we consider the following Fokker-Planck equation
\begin{equation}\label{Fokker-Planck}
d u(x, t) = \big( -  b \cdot \f{\p u(x, t)}{\p x} + \f{\delta^2 }{2} \cdot \f{\p^2 u(x, t)}{\p x^2} \big) dt,
\end{equation} 
which describes the time evolution of the probability density function of the velocity of a particle under the influence of drag forces and random forces. The parameters $b$ and $\delta$ in Eq. \eqref{Fokker-Planck} defines the drag force (drift) and the random force (diffusion), respectively.

To derive the state model, we discretize the above Fokker-Planck equation by using the following finite difference scheme with time step-size $\Delta t$ and spatial step-size $\Delta x$:
\begin{equation}\label{Dis:FP}
X^{(i)}_{n+1} = X^{(i)}_n - b \f{X^{(i+1)}_{n} - X^{(i)}_{n}}{\Delta x} \Delta t + \f{\delta^2}{2} \f{X^{(i+1)}_{n} - 2 X^{(i)}_{n} + X^{(i-1)}_{n}}{\Delta x^2} \Delta t + \sigma \xi^{(i)}_n,
\end{equation}
where $X^{(i)}_{n+1} \approx u(x_i, t_{n+1})$ is an approximation for solution $u$ at spatial point $x_i$ and time $t_{n+1},$
and we assume that the propagation of solution $u$ is perturbed by a time-space dependent noise $\sigma \xi^{(i)}_n$ with $\xi^{(i)}_n \sim N(0, 1)$ representing a standard Gaussian random variables for the given time instant $n$ corresponding to the spatial dimension discretization index $i = 1, 2, \cdots, d$.

The goal of the state-parameter estimation problem in this example is to use observational data to estimate the state of the approximated solution set $\{X_{n+1}^{(i)}\}_{i=1}^{d}$ at time instant $n+1$ as well as the unknown model parameters, i.e., $b$ and $\delta$. The observations $\{Y_{n+1}^{(i)}\}_{i=1}^{d}$ that we receive are defined by 
\begin{equation}\label{FP_obs}
Y_{n+1}^{(i)} = \max\{X^{(i)}_{n+1}, \lambda \} + \zeta^{(i)}_{n+1},
\end{equation}
where $\zeta^{(i)}_{n+1} \sim N(0, 0.02^2)$ is the observational noise at space point $x_i$ and time instant $n+1$, and $\lambda $ in the observation function is a cut-off threshold that reflects the sensitivity of the measurement detector -- only signal values greater than $\lambda $ can be received and measured by the detector. The choice of observation function in \eqref{FP_obs} allows us to examine the robustness of the United Filter when observations on the state model only provide partial information.

In our numerical experiments, we let the true unknown parameters be $b = 10$ and $\delta = 2$, and the initial guess values for these two parameters are $\bar{b}_0 = 2$ and $\bar{\delta}_0 = 10$. We let the initial condition of the Fokker-Planck equation be 
$$u(x, 0) = 10 \exp\big( - (x - 10)^2\big),$$
and the noise level is chosen as $\sigma = 0.5 \Delta t$ in \eqref{Dis:FP}. 
\begin{figure}[h!]
\begin{center}
\subfloat[State estimation: $X_{25}$]{\includegraphics[scale = 0.2]{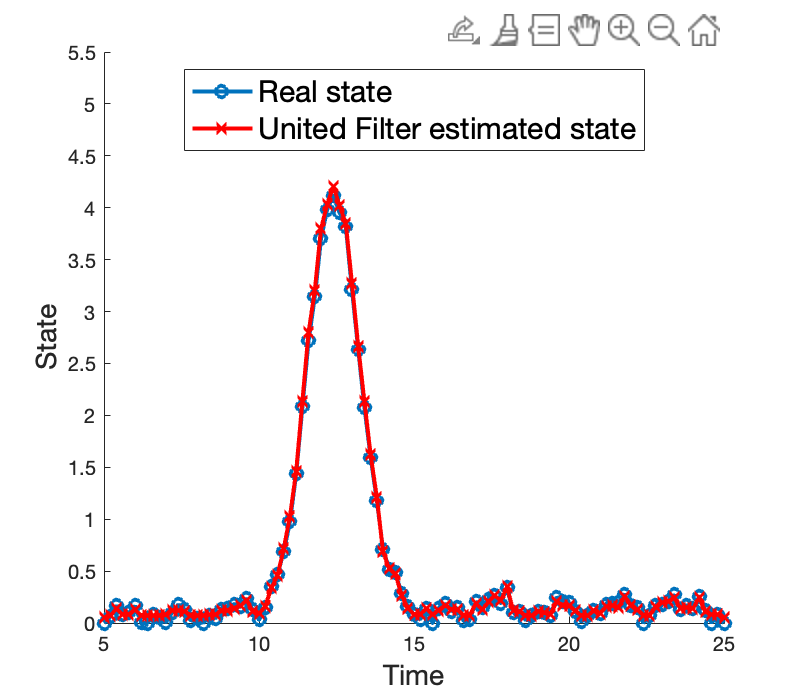} } \hspace{1em}
\subfloat[State estimation: $X_{50}$]{\includegraphics[scale = 0.2]{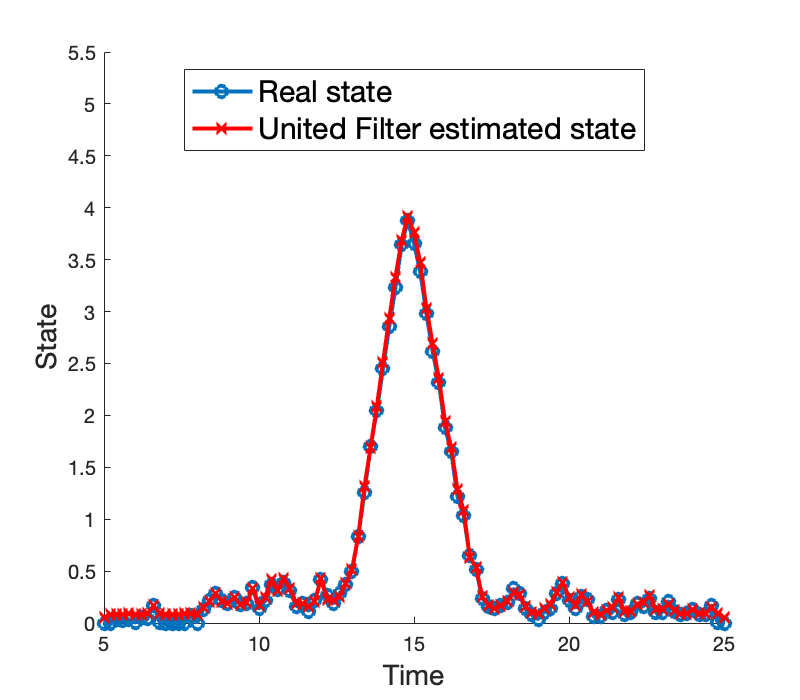} }   \\
\subfloat[State estimation: $X_{75}$]{\includegraphics[scale = 0.2]{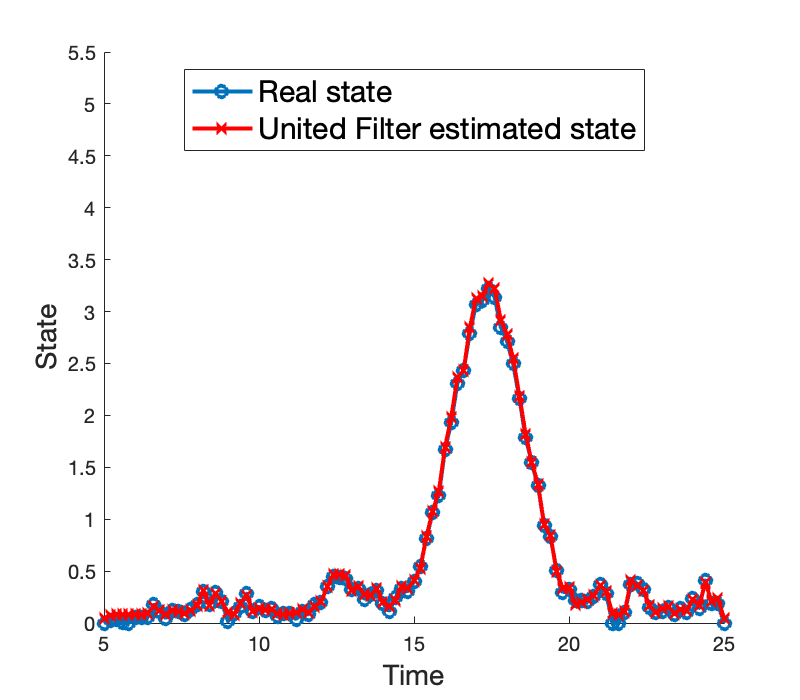} }   \hspace{1em}
\subfloat[State estimation: $X_{100}$]{\includegraphics[scale = 0.2]{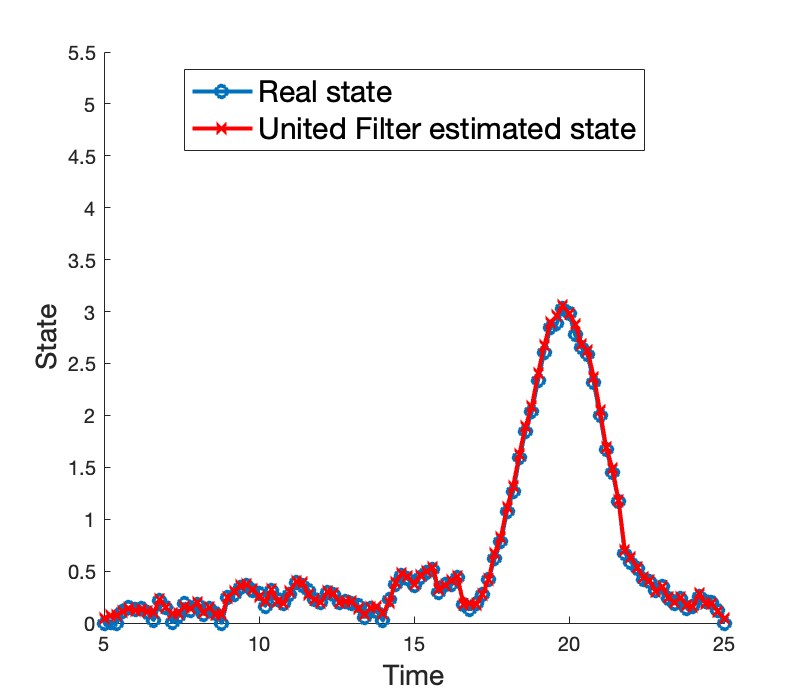} }   \\
\end{center} \vspace{-0.5em}
\caption{Example 2. State estimation of the United Filter with observation cut-off threshold $\lambda  = 0.1$. }\label{Ex1:FP_01_X}
\end{figure}
\begin{figure}[h!]
\begin{center}
\subfloat[Parameter estimation: $b$]{\includegraphics[scale = 0.25]{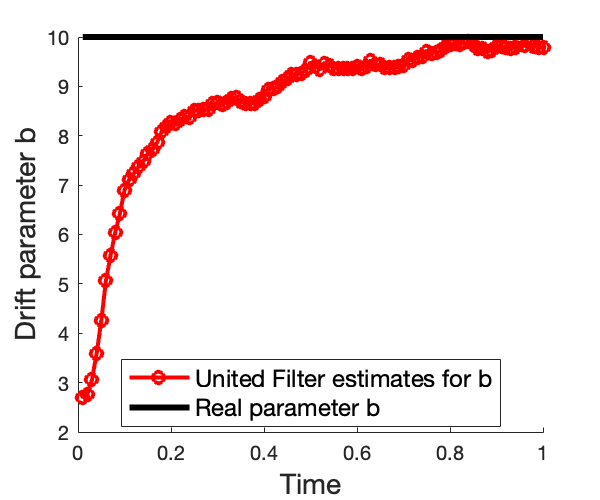} } \hspace{1em}
\subfloat[Parameter estimation: $\delta$]{\includegraphics[scale = 0.25]{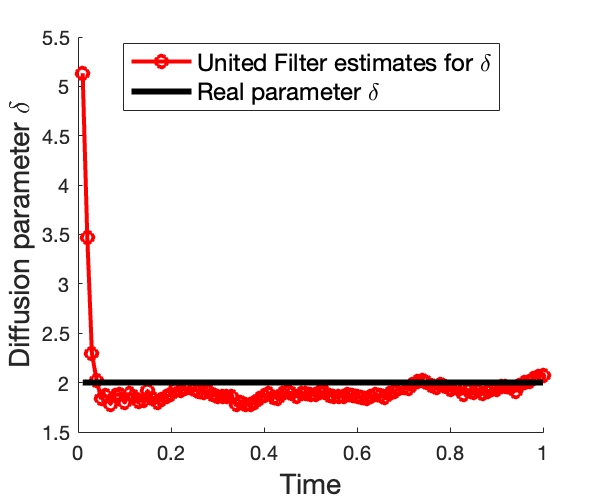} }  
\end{center} \vspace{-0.5em}
\caption{Example 2. Parameter estimation of the United Filter with observation cut-off threshold $\lambda  = 0.1$. }\label{Ex1:FP_01_para}
\end{figure}
To establish the finite difference scheme, we give a spatial partition to the $[0, 30]$ region with spatial step-size $\Delta x = 0.25$, i.e., the dimension of the problem is $d = 120$, and we estimate the state and parameter over the time interval $[0, 1]$ with time step-size $\Delta t = 0.01$.

In Figure \ref{Ex1:FP_01_X}, we show the performance of state estimation by using the United Filter at time instants $n = 25$, $n = 50$, $n = 75$, and $n = 100$, and we assume that the cut-off value to be $\lambda  = 0.1$. In Figure \ref{Ex1:FP_01_para}, we present the performance of parameter estimation over the $[0, 1]$ time interval. Note that the number of spatial partition is equivalent to the number of dimension in the state estimation problem. Therefore, the United Filter has solved a $120$-dimensional state estimation problem via the EnSF, and $d=120$ is a very high dimensional problem, which could challenge most existing optimal filtering methods. To implement the United Filter, we utilized $200$ diffusion model SDE (and reverse SDE) samples with $100$ time steps in the pseudo time interval $[0, 1]$, and we used $400$ particles to estimate the unknown parameters $b$ and $\delta$ in the direct filter component of the United Filter algorithm.
We can see from these figures that the EnSF component of the United Filter provided very accurate estimates for the state vector in the $120$-dimensional space, and the direct filter component for parameter estimation also provided reliable estimates for the unknown parameters.
\begin{figure}[h!]
\begin{center}
\subfloat[State estimation: $X_{25}$]{\includegraphics[scale = 0.2]{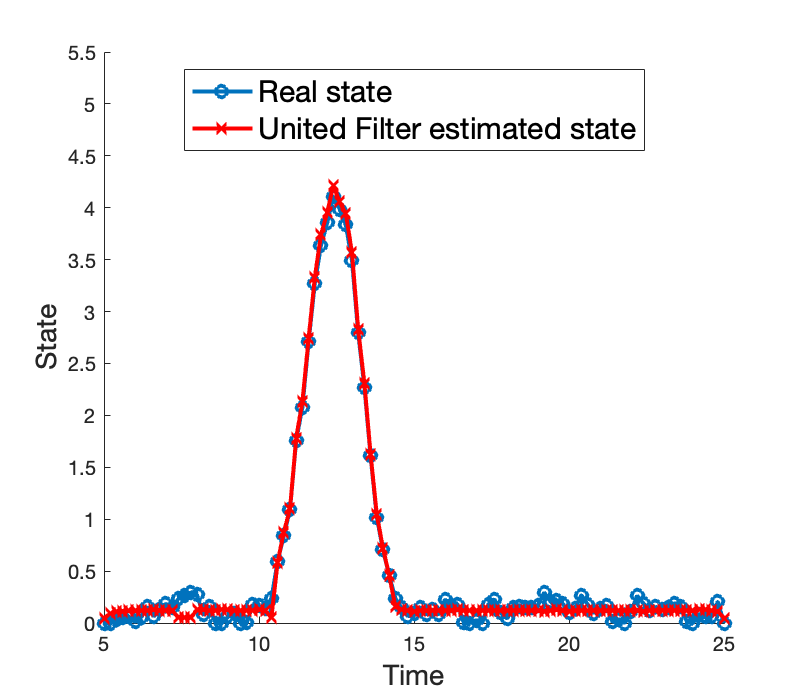} } \hspace{1em}
\subfloat[State estimation: $X_{50}$]{\includegraphics[scale = 0.2]{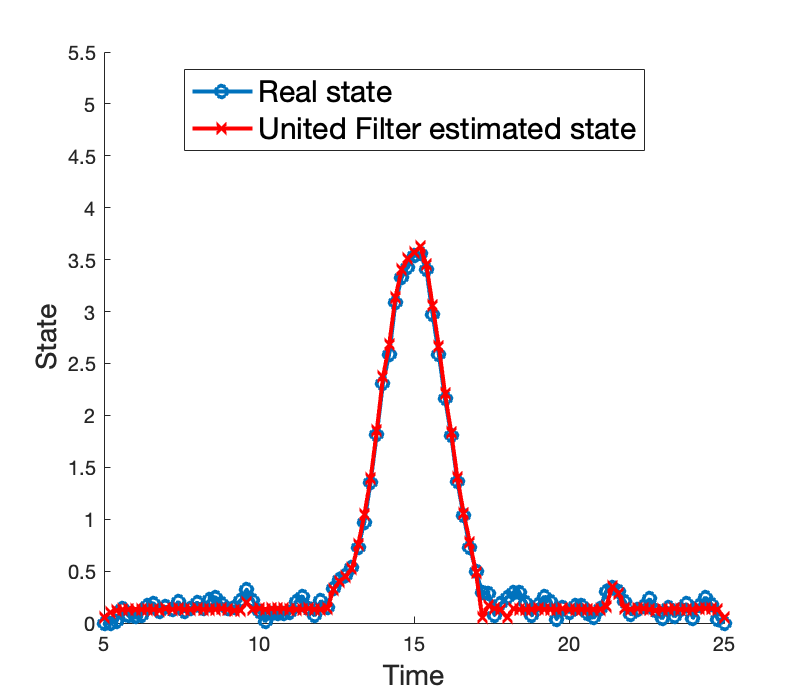} }   \\
\subfloat[State estimation: $X_{75}$]{\includegraphics[scale = 0.2]{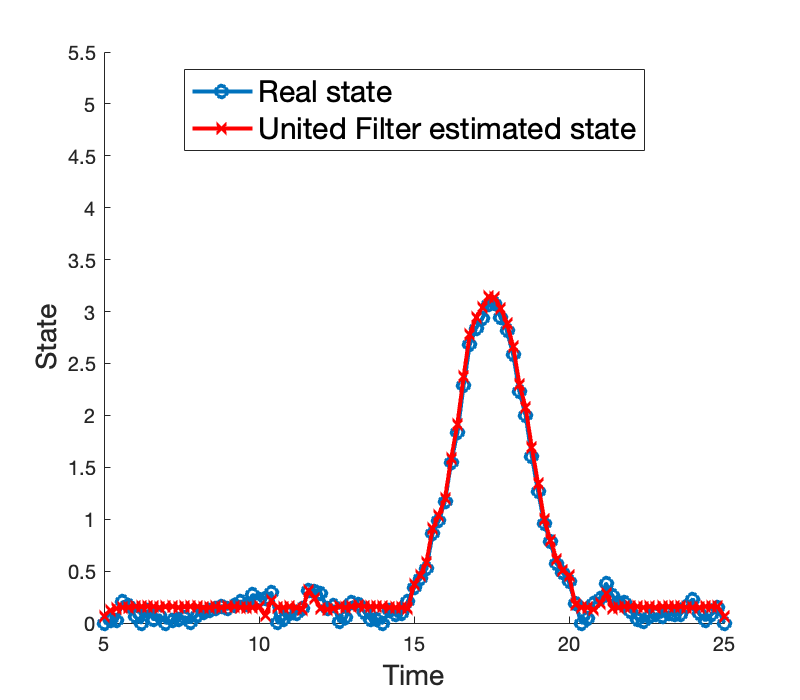} }   \hspace{1em}
\subfloat[State estimation: $X_{100}$]{\includegraphics[scale = 0.2]{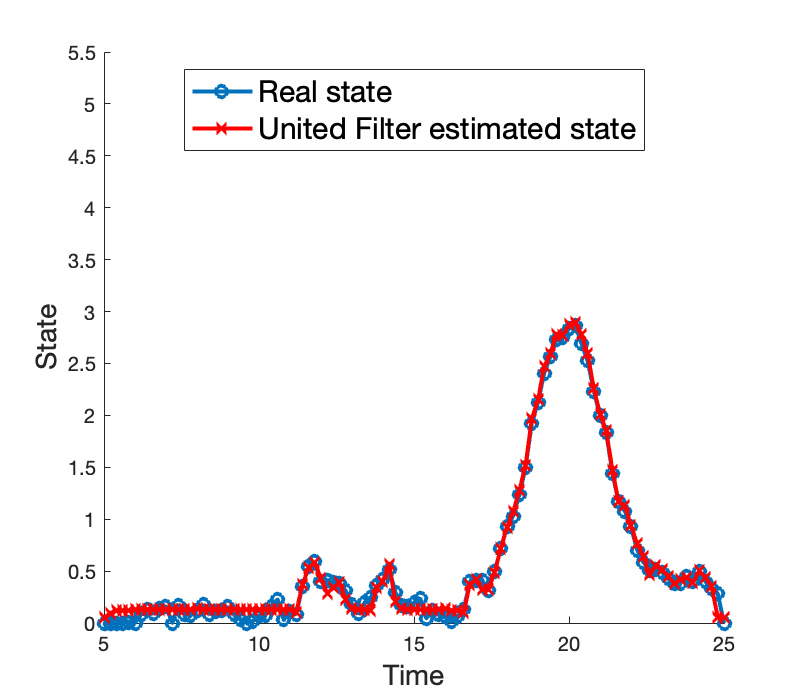} }   \\
\end{center} \vspace{-0.5em}
\caption{Example 2. State estimation of the United Filter with observation cut-off threshold $\lambda  = 0.3$. }\label{Ex1:FP_03_X}
\end{figure}
\begin{figure}[h!]
\begin{center}
\subfloat[Parameter estimation: $b$]{\includegraphics[scale = 0.25]{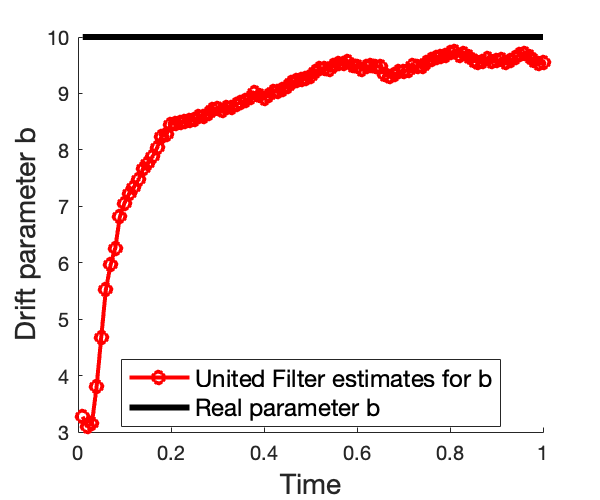} } \hspace{1em}
\subfloat[Parameter estimation: $\delta$]{\includegraphics[scale = 0.25]{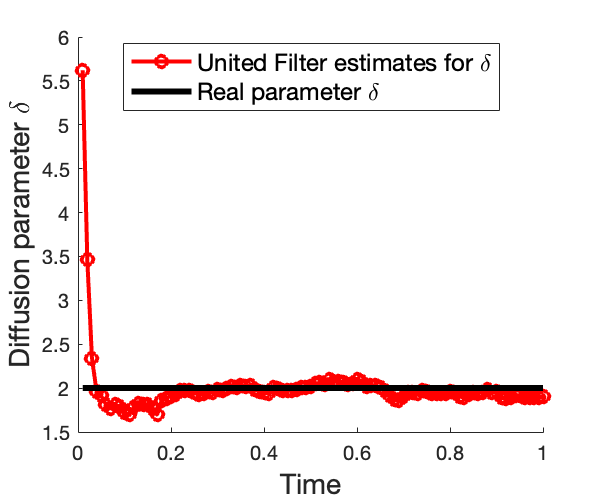} }  
\end{center} \vspace{-0.5em}
\caption{Example 2. Parameter estimation of the United Filter with observation cut-off threshold $\lambda  = 0.3$. }\label{Ex1:FP_03_para}
\end{figure}

To further demonstrate the robustness of our United Filter method, we increase the observation cut-off threshold to $\lambda  = 0.3$ and $\lambda  = 0.5$. 
\begin{figure}[h!]
\begin{center}
\subfloat[State estimation: $X_{25}$]{\includegraphics[scale = 0.2]{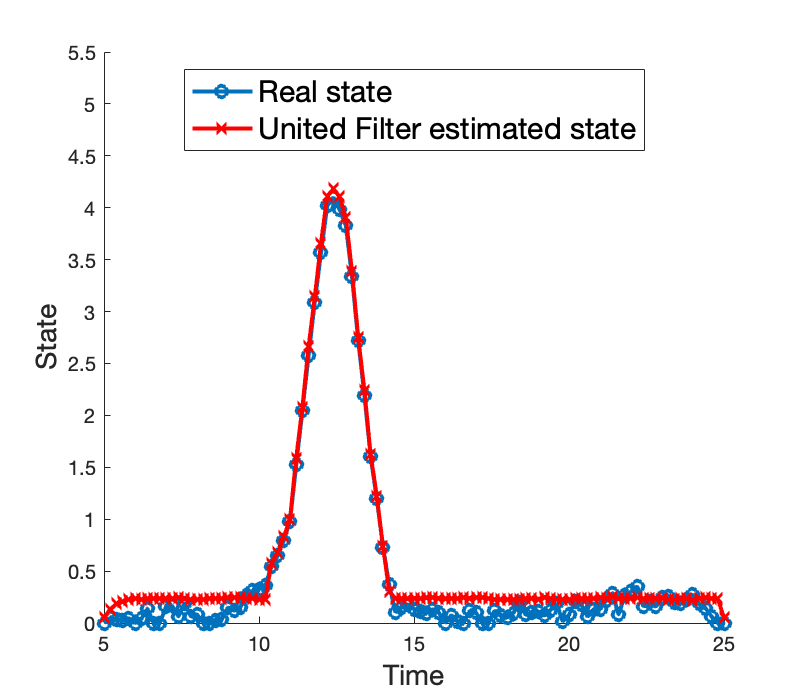} } \hspace{1em}
\subfloat[State estimation: $X_{50}$]{\includegraphics[scale = 0.2]{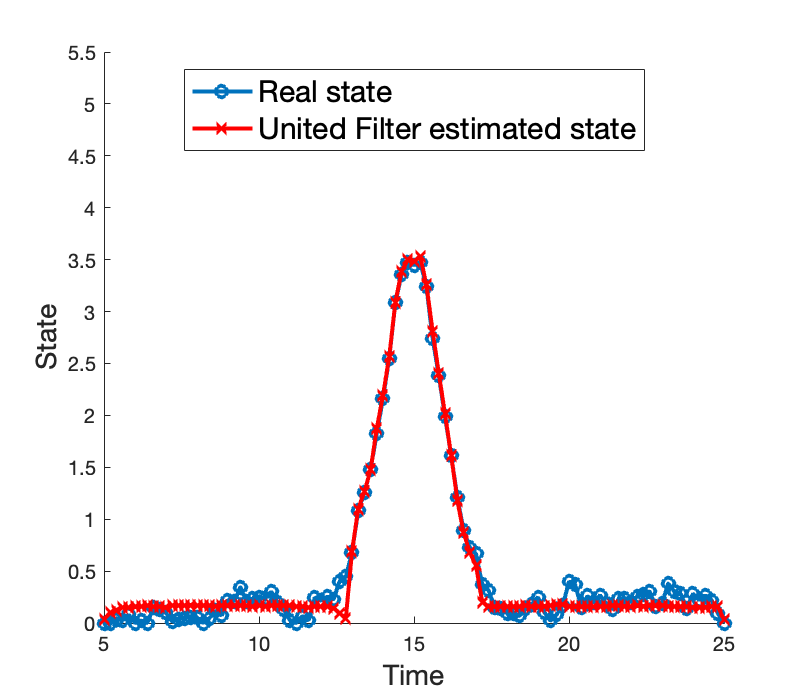} }   \\
\subfloat[State estimation: $X_{75}$]{\includegraphics[scale = 0.2]{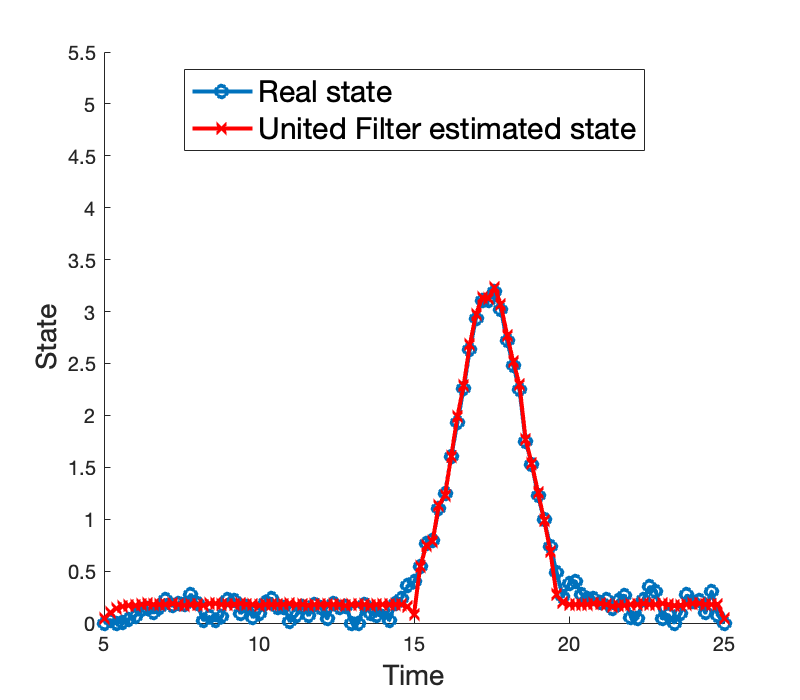} }   \hspace{1em}
\subfloat[State estimation: $X_{100}$]{\includegraphics[scale = 0.2]{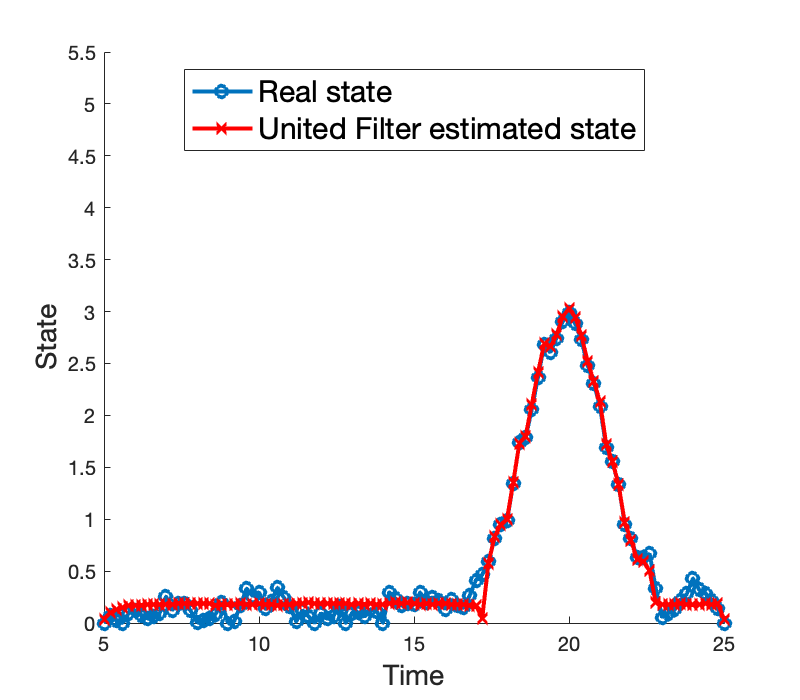} }   \\
\end{center} \vspace{-0.5em}
\caption{Example 2. State estimation of the United Filter with observation cut-off threshold $\lambda  = 0.5$. }\label{Ex1:FP_05_X}
\end{figure}
\begin{figure}[h!]
\begin{center}
\subfloat[Parameter estimation: $b$]{\includegraphics[scale = 0.25]{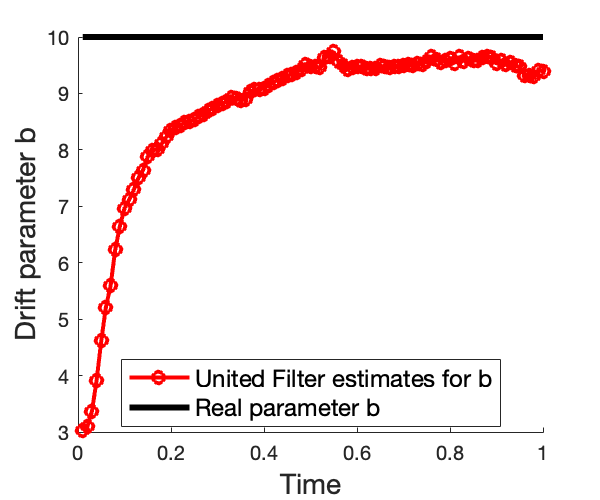} } \hspace{1em}
\subfloat[Parameter estimation: $\delta$]{\includegraphics[scale = 0.25]{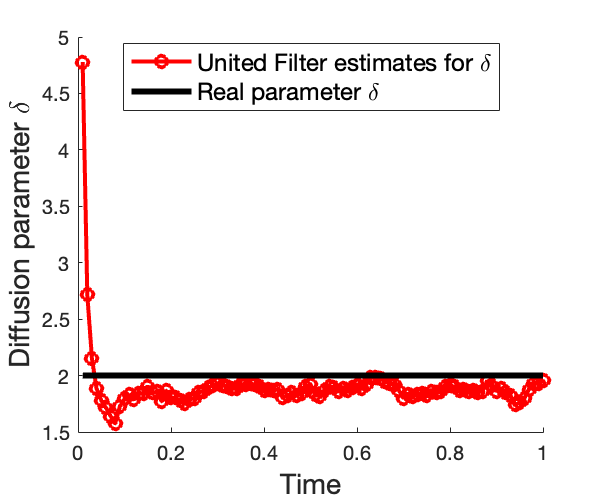} }  
\end{center} \vspace{-0.5em}
\caption{Example 2. Parameter estimation of the United Filter with observation cut-off threshold $\lambda  = 0.5$. }\label{Ex1:FP_05_para}
\end{figure}
In Figure \ref{Ex1:FP_03_X} and Figure \ref{Ex1:FP_03_para}, we show the state estimation performance and the parameter estimation performance obtained by the United Filter, respectively, in the case that $\lambda  = 0.3$. Similarly, in Figure \ref{Ex1:FP_05_X} and Figure \ref{Ex1:FP_05_para}, we show the state estimation performance and the parameter estimation performance in the case that $\lambda  = 0.5$.  From those figures, we can see that even with a less sensitive measurement detector that may overlook low-magnitude details, the state estimation provided by the EnSF still accurately captures the mode in the solution of the Fokker-Planck equation. Moreover, despite of higher observation cut-off thresholds, the United Filter can still provide good parameter estimation results.

\subsection{Example 3: Lorenz 96 model}
In this example, we solve the state-parameter estimation problem for the Lorenz 96 model, which is a well-known chaotic dynamical system that poses challenges for existing optimal filtering methods. When considering unknown parameters, this problem becomes even more difficult than standard state estimation. This is because predicting the model without true parameter values can lead to results that significantly deviate from the observational data, making high-dimensional state estimation much more challenging.

To show the advantageous performance of our method, we will compare the United Filter method with the Augmented Ensemble Kalman Filter (AugEnKF), which is the state-of-the-art approach for joint state-parameter estimation in high dimensional spaces.

The Lorenz 96 model that we consider is given as follows: 
\begin{equation}\label{Lorenz-state}
X^{(i)}_{n+1} = X^{(i)}_{n} +  \Big[ \lambda \big( X^{(i+1)}_{n} - X^{(i-2)}_{n}\big) X^{(i-1)}_{n} - \gamma X^{(i)}_{n} + F \Big] \Delta t  + \sigma \xi^{(i)}_{n}, \hspace{1em} i = 1, 2, \cdots, d,
\end{equation}
with $X^{(-1)}_{n} = X^{(d-1)}_{n}$, $X^{(0)}_{n} = X^{(d)}_{n}$ and $X^{(d+1)}_{n} = X^{(1)}_{n}$, $d \geq 4$, where $\lambda$, $\gamma$, $F$ are unknown parameters to be determined along with the state $\bm{X}_{n+1} = [X_{n+1}^{(1)}, X_{n+1}^{(2)}, \cdots, X_{n+1}^{(d)},]$, and $\Delta t$ is a time step-size. We assume that the Lorenz model \eqref{Lorenz-state} is perturbed by a $d$-dimensional standard Gaussian noise $\xi_n = \{\xi_n^{(i)}\}_{i=1}^d$ with coefficient $\sigma = 0.1 \sqrt{\Delta t}$. 

In this example, we set the true parameter values as $\lambda = 2$, $\gamma = 5$, and $F = 8$, while initializing estimates at $\bar{\lambda}_0 = 8$, $\bar{\gamma}_0 = 1$, and $\bar{F}_0 = 1$. The problem dimension is defined as $d=200$, and we perform state tracking and parameter estimation over 50 time steps, using a time step-size of $\Delta t = 0.02$.

To create a challenging joint state-parameter estimation problem, we consider that at each observation moment, we can only receive observations in $100$ randomly selected directions, and these observations are subject to Gaussian random noise with a standard deviation of $0.05 \cdot I_{100}$. Furthermore, to intensify the problem's complexity and introduce additional nonlinearity, we randomly select $10$ directions within the state observation vector and define the observations as the arctangent of $X$.
\begin{figure}[h!]
\begin{center}
\subfloat[State estimation: $X_{50}$]{\includegraphics[scale = 0.2]{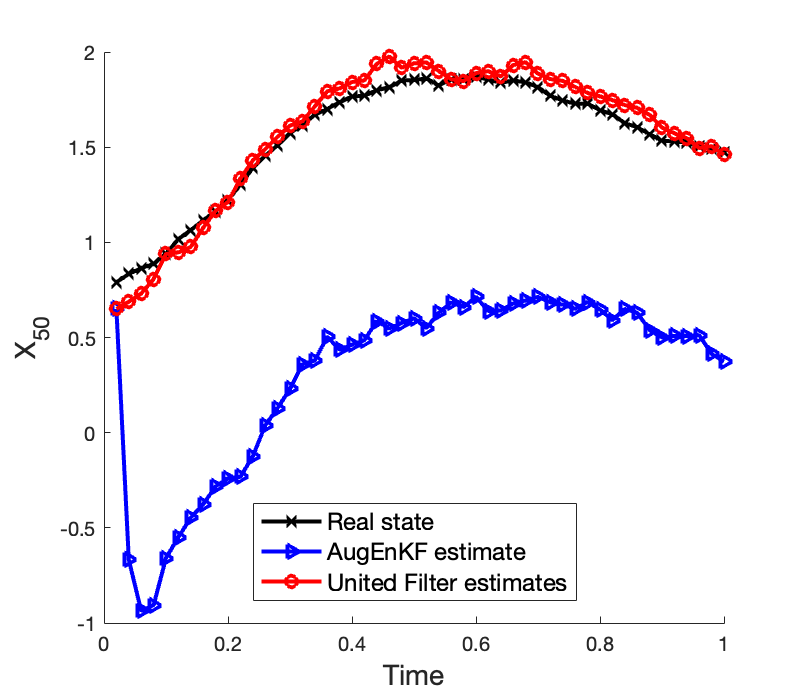} } \hspace{1em}
\subfloat[State estimation: $X_{100}$]{\includegraphics[scale = 0.2]{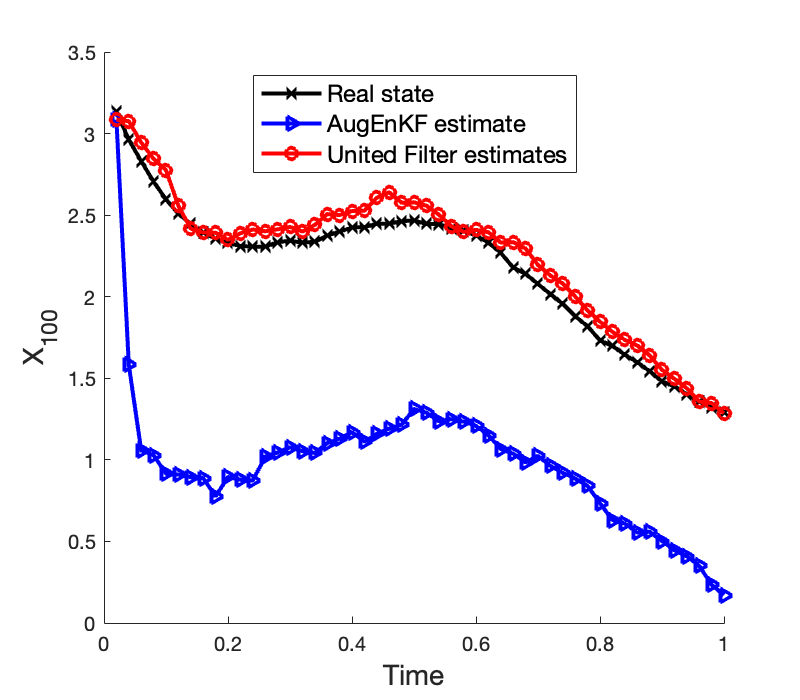} }   \\
\subfloat[State estimation: $X_{150}$]{\includegraphics[scale = 0.2]{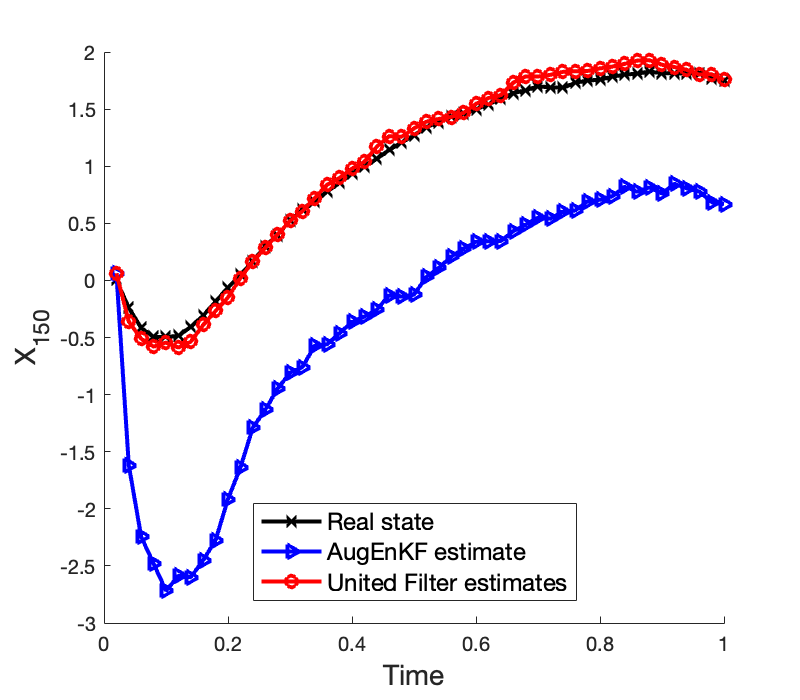} }   \hspace{1em}
\subfloat[State estimation: $X_{200}$]{\includegraphics[scale = 0.2]{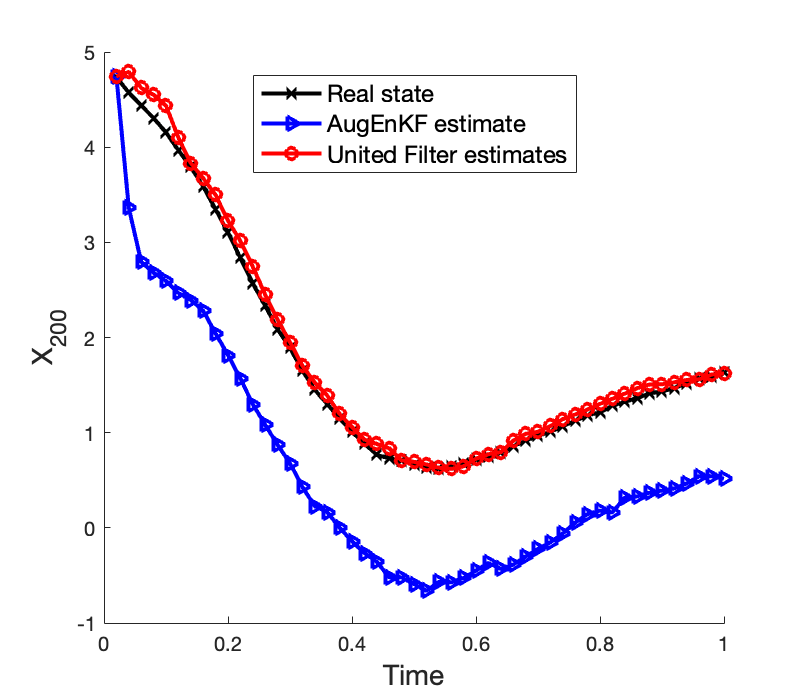} }   \\
\end{center} \vspace{-0.5em}
\caption{Example 3. Comparison of state estimation for $200$ dimensional Lorenz model. }\label{Ex3:Lorenz_200d_trajectory}
\end{figure}

In the forthcoming numerical experiments, we use $200$ diffusion SDE samples in the EnSF for conducting the state estimation procedure, and we use $1000$ particles to estimate the unknown parameters in the Lorenz 96 model. For the AugEnKF, we employ an ensemble of $1000$ realizations of the Kalman filter samples to carry out joint state-parameter estimation. It's important to note that working with $1000$ Kalman filter realizations is already considered a substantial number for solving a $203$-dimensional problem within the AugEnKF framework. 
In Figure \ref{Ex3:Lorenz_200d_trajectory}, we illustrate the state estimation comparison between the United Filter and the AugEnKF in dimensions $X_{50}$, $X_{100}$, $X_{150}$, and $X_{200}$. This figure clearly demonstrates the superior performance of the United Filter over the AugEnKF.
%
\begin{figure}[h!]
\begin{center}
\includegraphics[scale = 0.25]{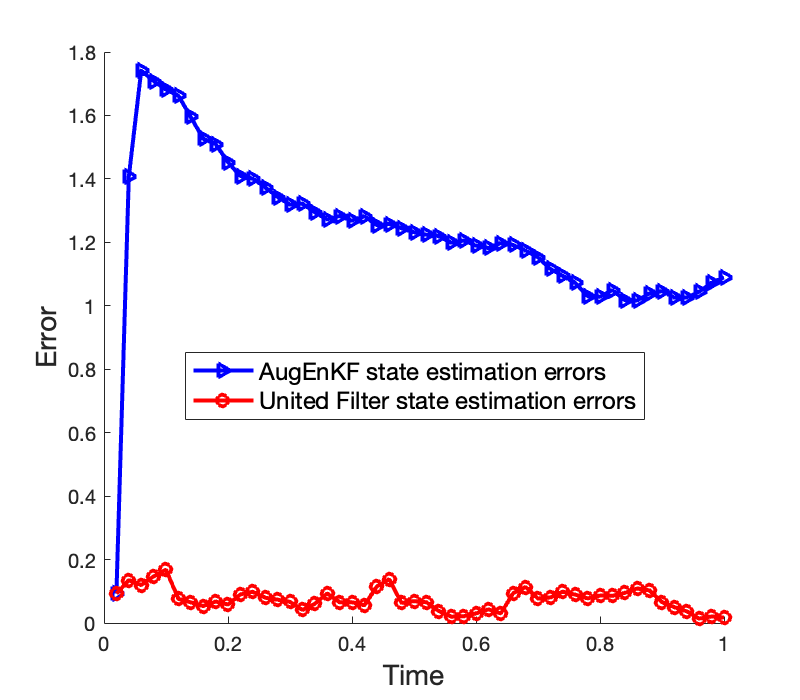} 
\end{center} \vspace{-0.5em}
\caption{Example 3. Comparison of average errors of state estimation in $200$ dimensions.}\label{Ex3:Lorenz_state}
\end{figure}
To assess the overall accuracy of state estimation, we calculated average errors across all $200$ directions, and the comparison result is depicted in Figure \ref{Ex3:Lorenz_state}. This result confirms the United Filter's superiority in state estimation.

\begin{figure}[h!]
\begin{center}
\subfloat[Parameter estimation for $\lambda$ ]{\includegraphics[scale = 0.18]{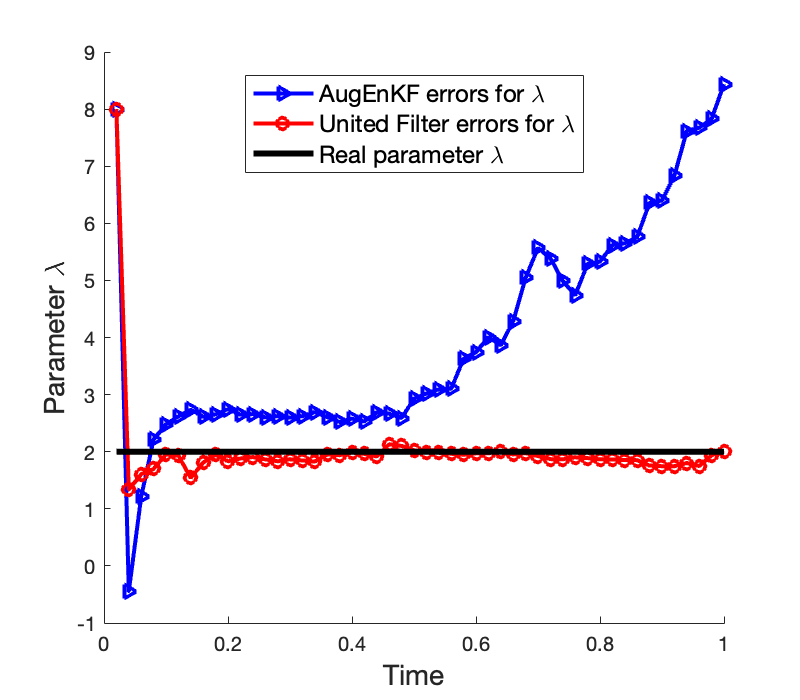} } 
\subfloat[Parameter estimation for $\gamma$]{\includegraphics[scale = 0.18]{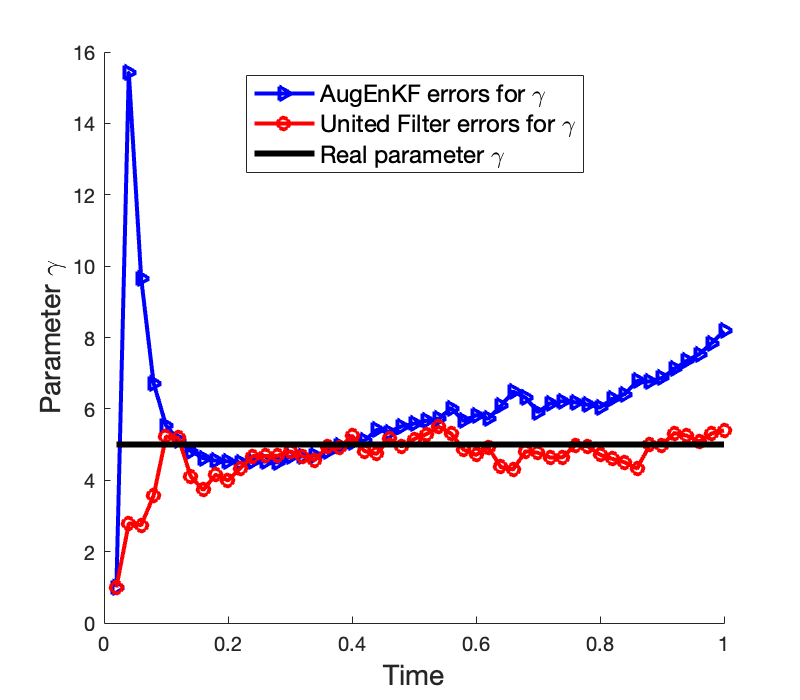} }   
\subfloat[Parameter estimation for $F$]{\includegraphics[scale = 0.18]{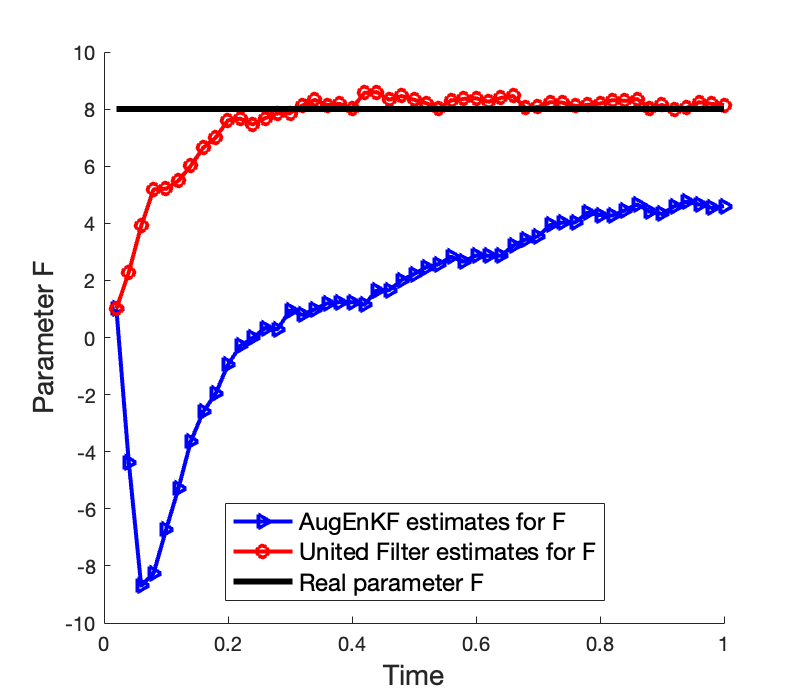} }
\end{center} \vspace{-0.5em}
\caption{Example 3. Comparison of parameter estimation for the $200$ dimensional Lorenz model. }\label{Ex3:Lorenz_parameter}
\end{figure}

In Figure \ref{Ex3:Lorenz_parameter}, we compare the performance of parameter estimation between the United Filter and the AugEnKF, where subplots (a), (b), and (c) show parameter estimation results for $\lambda$, $\gamma$, and $F$, respectively.
We can see from this result that the United Filter also outperforms the AugEnKF in parameter estimation. 

\begin{figure}[h!]
\begin{center}
\subfloat[RMSEs for estimating $\lambda$]{\includegraphics[scale = 0.18]{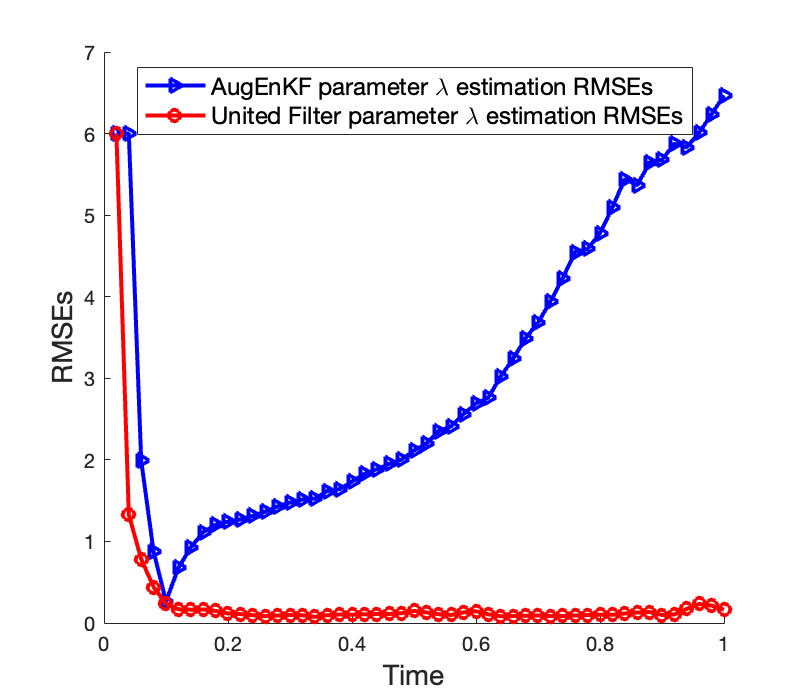} } 
\subfloat[RMSEs for estimating $\gamma$]{\includegraphics[scale = 0.18]{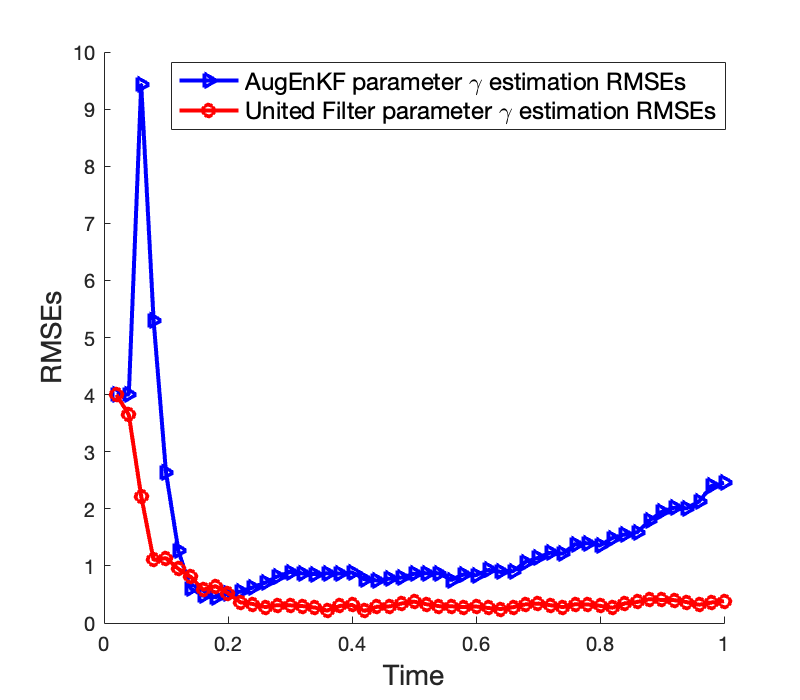} }   
\subfloat[RMSEs for estimating $F$]{\includegraphics[scale = 0.18]{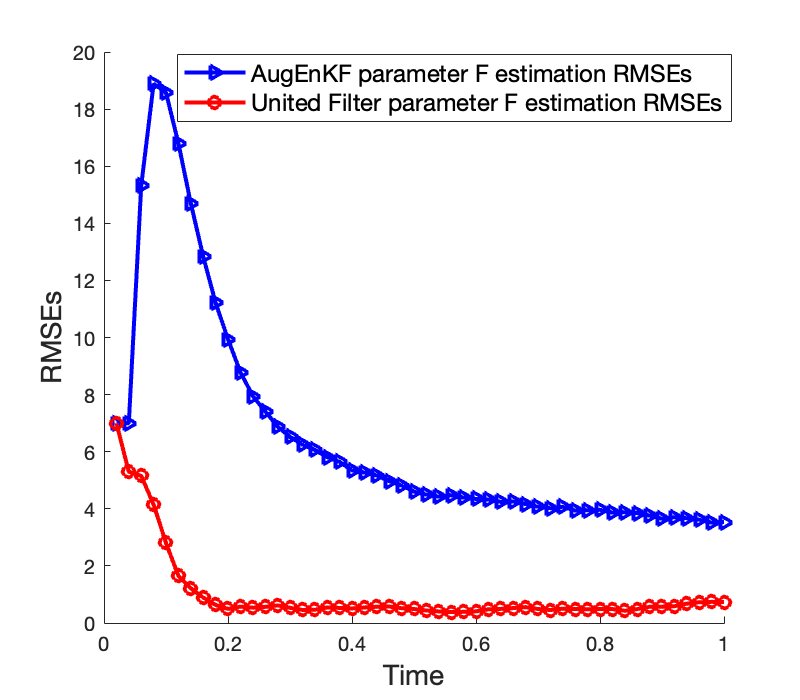} }
\end{center} \vspace{-0.5em}
\caption{Example 3. Comparison of RMSEs in parameter estimation. }\label{Ex3:Lorenz_parameter_RMSE}
\end{figure}

Note that in the state-parameter estimation problem, the accuracy of state estimation and parameter estimation are interrelated. Poor parameter estimation can adversely affect state estimation accuracy, while imprecise state estimation can similarly lead to inaccurate parameter estimation.
\begin{figure}[h!]
\begin{center}
\includegraphics[scale = 0.25]{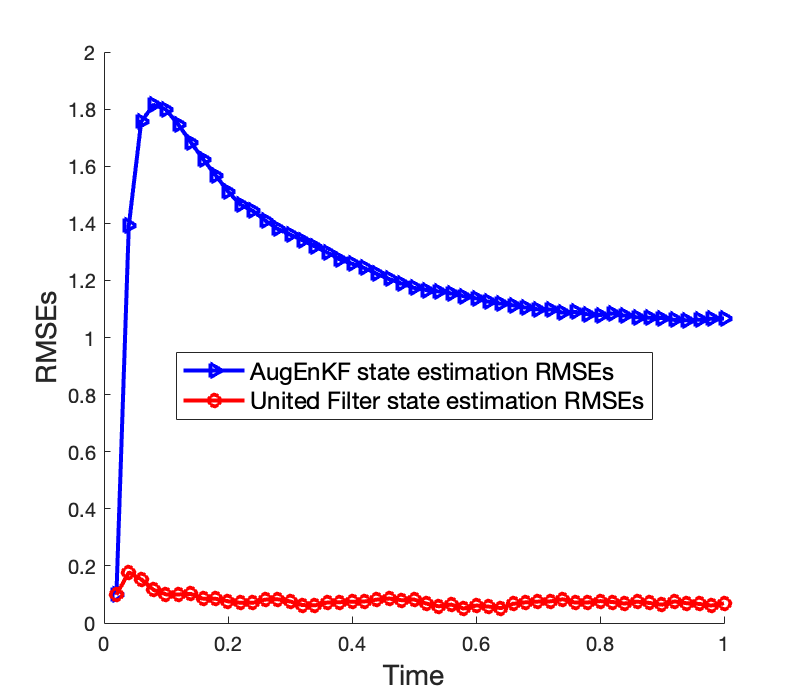} 
\end{center} \vspace{-0.5em}
\caption{Example 3. Comparison of RMSEs in state estimation.}\label{Ex3:Lorenz_state_RMSE}
\end{figure}
We validate the results from the above numerical experiments by conducting 20 repetitions and calculating the Root Mean Square Errors (RMSEs). In Figure \ref{Ex3:Lorenz_state_RMSE}, we compare the RMSEs in state estimation between the United Filter and the AugEnKF, and in Figure \ref{Ex3:Lorenz_parameter_RMSE} we compare the RMSEs in parameter estimation. 
As we can observe in Figure \ref{Ex3:Lorenz_state_RMSE} and Figure \ref{Ex3:Lorenz_parameter_RMSE}, the United Filter constantly outperforms the AugEnKF in both state estimation and parameter estimation.
\begin{figure}[h!]
\begin{center}
\subfloat[AugEnKF estimation trajectories for $\lambda$]{\includegraphics[scale = 0.24]{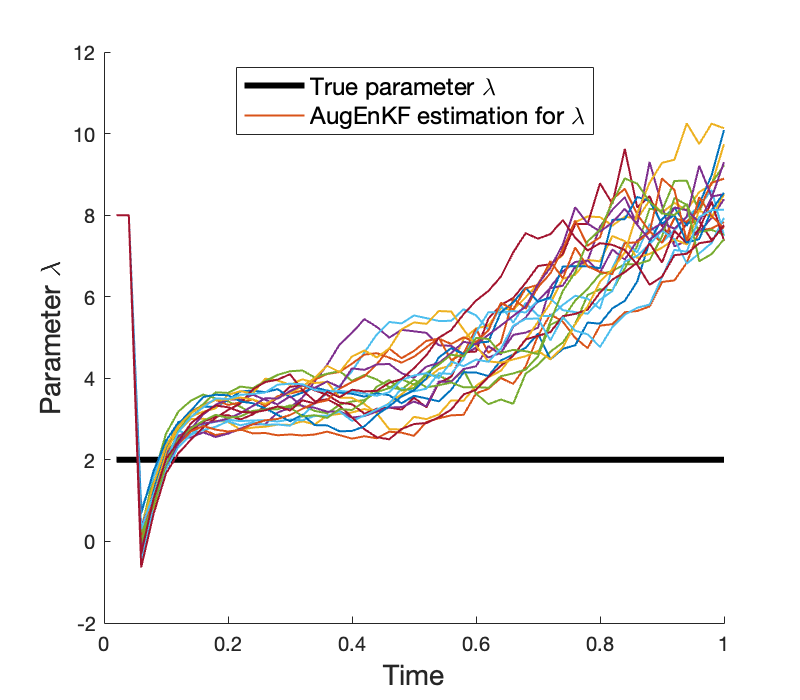} } \hspace{1em}
\subfloat[United Filter estimation trajectories for $\lambda$]{\includegraphics[scale = 0.24]{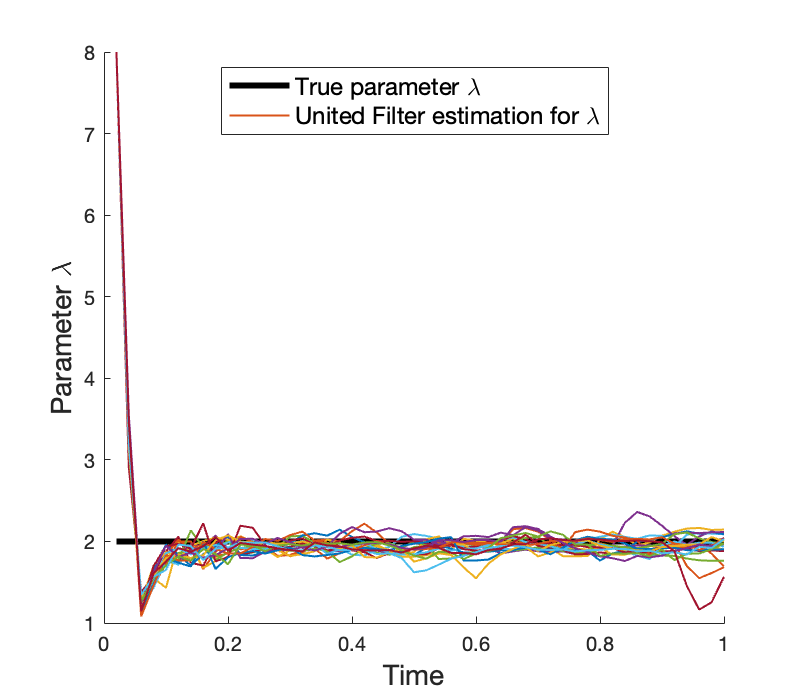} }   \\
\end{center} \vspace{-0.5em}
\caption{Example 3. Comparison of parameter estimation in $20$ repeated experiments: $\lambda$. }\label{Ex3:Lorenz_lambda_trajectory}
\end{figure}
\begin{figure}[h!]
\begin{center}
\subfloat[AugEnKF estimation trajectories for $\gamma$]{\includegraphics[scale = 0.24]{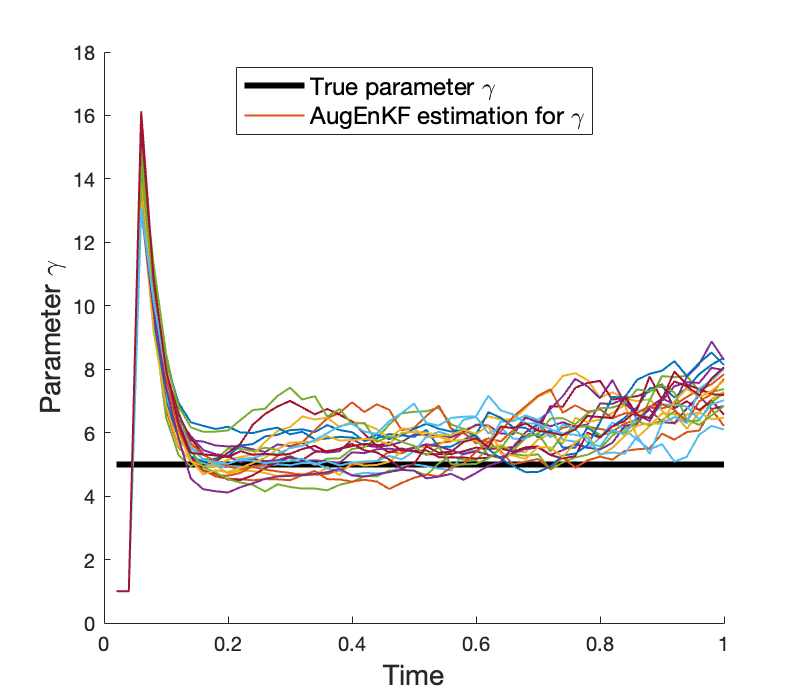} } \hspace{1em}
\subfloat[United Filter estimation trajectories for $\gamma$]{\includegraphics[scale = 0.24]{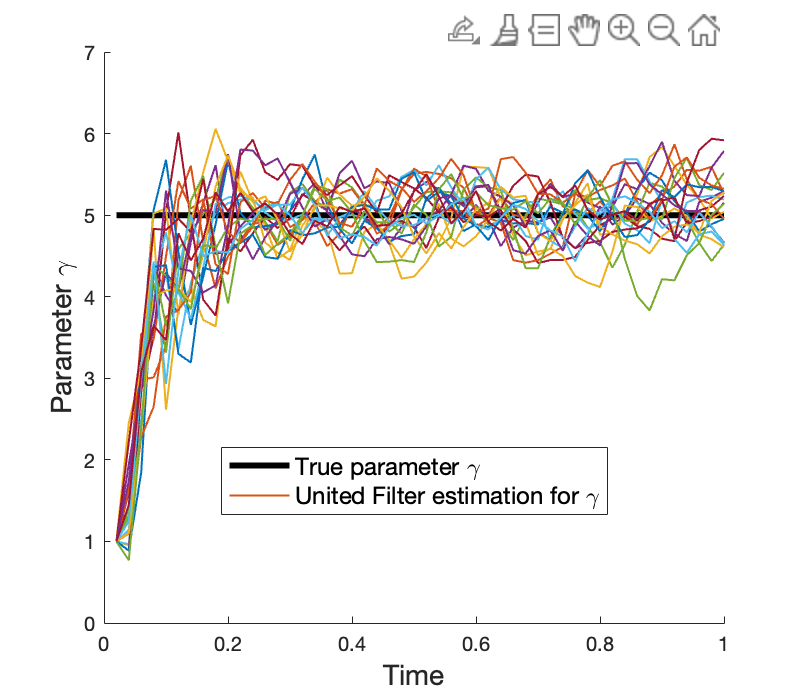} }   \\
\end{center} \vspace{-0.5em}
\caption{Example 3. Comparison of parameter estimation in $20$ repeated experiments: $\gamma$. }\label{Ex3:Lorenz_gamma_trajectory}
\end{figure}
\begin{figure}[h!]
\begin{center}
\subfloat[AugEnKF estimation trajectories for $F$]{\includegraphics[scale = 0.24]{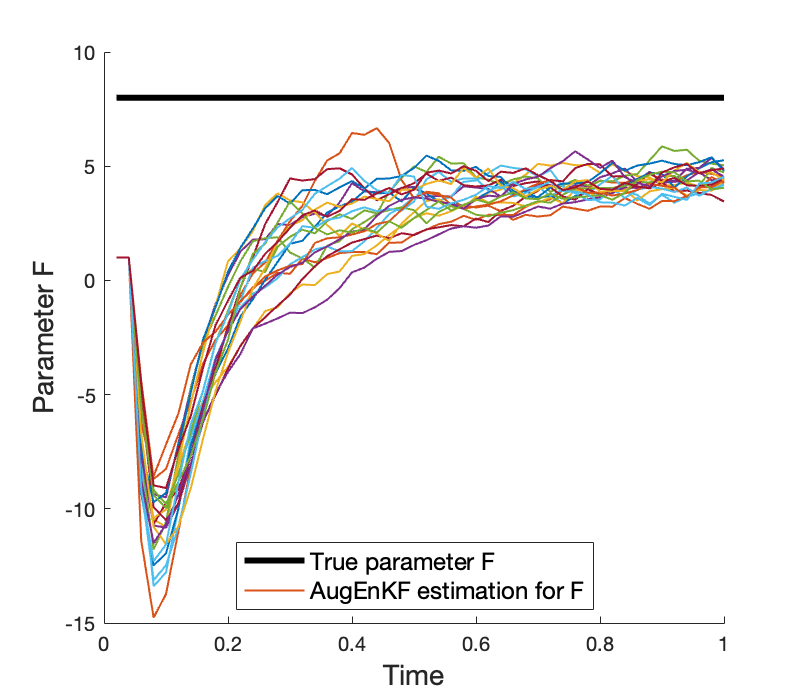} } \hspace{1em}
\subfloat[United Filter estimation trajectories for $F$]{\includegraphics[scale = 0.24]{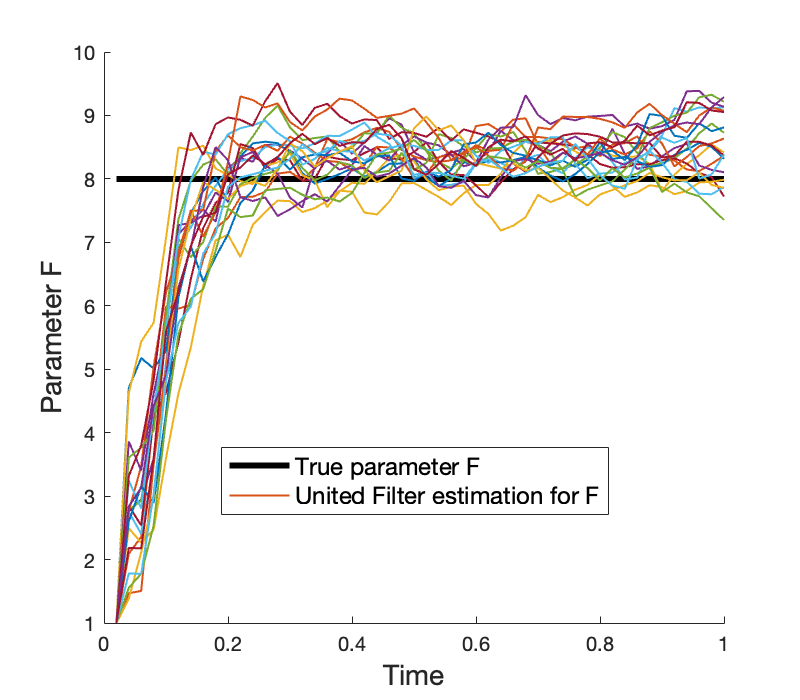} }   \\
\end{center} \vspace{-0.5em}
\caption{Example 3. Comparison of parameter estimation in $20$ repeated experiments: $F$. }\label{Ex3:Lorenz_F_trajectory}
\end{figure}

To provide a more comprehensive view of parameter estimation performance across $20$ repeated tests, we display all $20$ trajectories of parameter estimation results obtained by both the United Filter and the AugEnKF for $\lambda$, $\gamma$ and $F$ in Figures \ref{Ex3:Lorenz_lambda_trajectory},  \ref{Ex3:Lorenz_gamma_trajectory} and \ref{Ex3:Lorenz_F_trajectory}, respectively. We can see from these parameter estimation trajectories that the United Filter constantly generates very reliable estimates for all three parameters, with fluctuations well within a reasonable range.  In contrast, the AugEnKF fails to provide sufficiently accurate estimates for the unknown parameters.

\section*{Acknowledgement.}
This material is based upon work supported by the U.S. Department of Energy, Office of Science, Office of Advanced Scientific Computing Research, Applied Mathematics program under the contract ERKJ387 at the Oak Ridge National Laboratory, which is operated by UT-Battelle, LLC, for the U.S. Department of Energy under Contract DE-AC05-00OR22725. The First author (FB) would also like to acknowledge the support from U.S. National Science Foundation through
project DMS-2142672 and the support from the U.S. Department of Energy, Office of Science, Office of Advanced Scientific Computing Research, Applied Mathematics program under Grant DE-SC0022297.

 \bibliographystyle{elsarticle-num}

\end{document}